\documentclass[10pt]{amsart}
\usepackage{amscd}

\DeclareFontEncoding{OT2}{}{} 
  \newcommand{\textcyr}[1]{%
    {\fontencoding{OT2}\fontfamily{wncyr}\fontseries{m}\fontshape{n}%
     \selectfont #1}}
\newcommand{\Sha}{{\mbox{\textcyr{Sh}}}}

\address{1126 Burnside Hall \\ Department of Mathematics and Statistics \\
McGill University \\ 805 Sherbrooke West \\ Montreal, QC, Canada H3A 2K6}
\email{clark@math.mcgill.ca}

\title{The Period-Index Problem in WC-groups II: abelian varieties}
\author{Pete L. Clark}
\begin{document}
\newtheorem{lemma}{Lemma}
\newtheorem{prop}[lemma]{Proposition}
\newtheorem{cor}[lemma]{Corollary}
\newtheorem{thm}[lemma]{Theorem}
\newtheorem{conj}[lemma]{Conjecture}
\newtheorem{quest}{Question}
\newtheorem{ques}[lemma]{Question}
\newtheorem{prob}{Problem}
\maketitle

\newcommand{\F}{\ensuremath{\mathbb F}}
\newcommand{\Fp}{\ensuremath{\F_p}}
\newcommand{\Fpbar}{\overline{\Fp}}
\newcommand{\Fq}{\ensuremath{\F_q}}
\newcommand{\N}{\ensuremath{\mathbb N}}
\newcommand{\Q}{\ensuremath{\mathbb Q}}
\newcommand{\R}{\ensuremath{\mathbb R}}
\newcommand{\Z}{\ensuremath{\mathbb Z}}
\newcommand{\SSS}{\ensuremath{\mathcal{S}}}
\newcommand{\Rn}{\ensuremath{\mathbb R^n}}
\newcommand{\Ri}{\ensuremath{\R^\infty}}
\newcommand{\C}{\ensuremath{\mathbb C}}
\newcommand{\Cn}{\ensuremath{\mathbb C^n}}
\newcommand{\Ci}{\ensuremath{\C^\infty}}
\newcommand{\U}{\ensuremath{\mathcal U}}
\newcommand{\gn}{\ensuremath{\gamma^n}}
\newcommand{\ra}{\ensuremath{\rightarrow}}
\newcommand{\fhat}{\ensuremath{\hat{f}}}
\newcommand{\ghat}{\ensuremath{\hat{g}}}
\newcommand{\hhat}{\ensuremath{\hat{h}}}
\newcommand{\covui}{\ensuremath{\{U_i\}}}
\newcommand{\covvi}{\ensuremath{\{V_i\}}}
\newcommand{\covwi}{\ensuremath{\{W_i\}}}
\newcommand{\Gt}{\ensuremath{\tilde{G}}}
\newcommand{\gt}{\ensuremath{\tilde{\gamma}}}
\newcommand{\Gtn}{\ensuremath{\tilde{G_n}}}
\newcommand{\gtn}{\ensuremath{\tilde{\gamma_n}}}
\newcommand{\gnt}{\ensuremath{\gtn}}
\newcommand{\Gnt}{\ensuremath{\Gtn}}
\newcommand{\Cpi}{\ensuremath{\C P^\infty}}
\newcommand{\Cpn}{\ensuremath{\C P^n}}
\newcommand{\lla}{\ensuremath{\longleftarrow}}
\newcommand{\lra}{\ensuremath{\longrightarrow}}
\newcommand{\Rno}{\ensuremath{\Rn_0}}
\newcommand{\dlra}{\ensuremath{\stackrel{\delta}{\lra}}}
\newcommand{\pii}{\ensuremath{\pi^{-1}}}
\newcommand{\la}{\ensuremath{\leftarrow}}
\newcommand{\gonem}{\ensuremath{\gamma_1^m}}
\newcommand{\gtwon}{\ensuremath{\gamma_2^n}}
\newcommand{\omegabar}{\ensuremath{\overline{\omega}}}
\newcommand{\dlim}{\underset{\lra}{\lim}}
\newcommand{\ilim}{\operatorname{\underset{\lla}{\lim}}}
\newcommand{\Hom}{\operatorname{Hom}}
\newcommand{\Ext}{\operatorname{Ext}}
\newcommand{\Part}{\operatorname{Part}}
\newcommand{\Ker}{\operatorname{Ker}}
\newcommand{\im}{\operatorname{im}}
\newcommand{\ord}{\operatorname{ord}}
\newcommand{\unr}{\operatorname{unr}}
\newcommand{\B}{\ensuremath{\mathcal B}}
\newcommand{\Ocr}{\ensuremath{\Omega_*}}
\newcommand{\Rcr}{\ensuremath{\Ocr \otimes \Q}}
\newcommand{\Cptwok}{\ensuremath{\C P^{2k}}}
\newcommand{\CC}{\ensuremath{\mathcal C}}
\newcommand{\gtkp}{\ensuremath{\tilde{\gamma^k_p}}}
\newcommand{\gtkn}{\ensuremath{\tilde{\gamma^k_m}}}
\newcommand{\PP}{\ensuremath{\mathcal P}}
\newcommand{\QQ}{\ensuremath{\mathcal Q}}
\newcommand{\I}{\ensuremath{\mathcal I}}
\newcommand{\sbar}{\ensuremath{\overline{s}}}
\newcommand{\Kn}{\ensuremath{\overline{K_n}^\times}}
\newcommand{\tame}{\operatorname{tame}}
\newcommand{\Qpt}{\ensuremath{\Q_p^{\tame}}}
\newcommand{\Qpu}{\ensuremath{\Q_p^{\unr}}}
\newcommand{\scrT}{\ensuremath{\mathfrak{T}}}
\newcommand{\That}{\ensuremath{\hat{\mathfrak{T}}}}
\newcommand{\Gal}{\operatorname{Gal}}
\newcommand{\Aut}{\operatorname{Aut}}
\newcommand{\tors}{\operatorname{tors}}
\newcommand{\Zhat}{\hat{\Z}}
\newcommand{\linf}{\ensuremath{l_\infty}}
\newcommand{\Lie}{\operatorname{Lie}}
\newcommand{\GL}{\operatorname{GL}}
\newcommand{\EEnd}{\operatorname{End}}
\newcommand{\aone}{\ensuremath{(a_1,\ldots,a_k)}}
\newcommand{\raone}{\ensuremath{r(a_1,\ldots,a_k,N)}}
\newcommand{\rtwoplus}{\ensuremath{\R^{2  +}}}
\newcommand{\rkplus}{\ensuremath{\R^{k +}}}
\newcommand{\length}{\operatorname{length}}
\newcommand{\Vol}{\operatorname{Vol}}
\newcommand{\cross}{\operatorname{cross}}
\newcommand{\GoN}{\Gamma_0( N)}
\newcommand{\DD}{\cal D}
\newcommand{\GAG}{\Gamma \alpha \Gamma}
\newcommand{\GBG}{\Gamma \beta \Gamma}
\newcommand{\HGD}{H(\Gamma,\Delta)}
\newcommand{\Ga}{\Gamma^\alpha}
\newcommand{\Div}{\operatorname{Div}}
\newcommand{\Divo}{\Div_0}
\newcommand{\Hstar}{\cal{H}^*}
\newcommand{\txon}{\tilde{X}_0(N)}

\newcommand{\abcd}{\left[ \begin{array}{cc}
a & b \\
c & d
\end{array} \right]}
\newcommand{\w}{\omega}

\newcommand{\Qpi}{\ensuremath{\Q(\pi)}}
\newcommand{\Qpin}{\Q(\pi^n)}
\newcommand{\End}{\operatorname{End}}
\newcommand{\pibar}{\overline{\pi}}
\newcommand{\pbar}{\overline{p}}
\newcommand{\lcm}{\operatorname{lcm}}
\newcommand{\trace}{\operatorname{trace}}
\newcommand{\OKv}{\mathcal{O}_{K_v}}
\newcommand{\OO}{\mathcal{O}}
\newcommand{\Pic}{\operatorname{Pic}}
\newcommand{\FPic}{\mathbf{Pic}}
\newcommand{\Gm}{\mathbb G_m}
\renewcommand{\Div}{\operatorname{Div}}
\newcommand{\gk}{\mathfrak{g}_k}
\newcommand{\normresn}{\langle \ , \ \rangle_n} 
\newcommand{\normresp}{\langle \ , \ \rangle_p}
\newcommand{\normcc}{\langle C_1, C_2 \rangle_p}
\newcommand{\Proof}{\begin{bf} Proof: \end{bf}}
\newcommand{\sep}{\operatorname{sep}}
\newcommand{\solv}{\operatorname{solv}}
\newcommand{\Alb}{\operatorname{Alb}}
\newcommand{\FAlb}{\mathbf{Alb}}
\newcommand{\gK}{\mathfrak{g}_k}
\newcommand{\kbar}{\overline{k}}
\newcommand{\Kbar}{\overline{k}}
\newcommand{\ggg}{\mathfrak{g}}
\newcommand{\jku}{\operatorname{YU}}
\newcommand{\car}{\operatorname{char}}
\begin{abstract}
We study the relationship between the period and the 
index of a principal homogeneous space over an abelian variety, 
obtaining results which generalize work of Cassels and Lichtenbaum on
curves of genus one.  In addition,
we show that the $p$-torsion in the Shafarevich-Tate group of a
fixed abelian variety over a number field $k$ grows arbitarily large
when considered over field extensions $l/k$ of bounded degree.  Essential use 
is made of an abelian variety version of O'Neil's period-index obstruction.
\end{abstract}
\noindent

\section{Introduction and Statement of Results}
Let $A/k$ be an abelian variety and $V/k$ a principal homogeneous space for $A$,\footnote{We should 
write $(V,\mu)$, where $\mu: A \times V \ra V$ is the structure map for $V$.  For the sake of
notational simplicity we neglect the $\mu$\ldots} so
$V$ represents an element of the Weil-Ch\^atelet group $H^1(k,A)$.  In this paper we are concerned
with two numerical invariants associated to $V$, each of which can be viewed as measuring the failure of $V$ to
 have $k$-rational points.
The \textbf{period} of $V$ is the order of $V$ in $H^1(k,A)$ (a torsion group), whereas the \textbf{index} of
$V$ is the greatest common divisor of all degrees of field extensions $[l:k]$ such that $V(l) \neq \emptyset$.
\footnote{\ldots but note that
the period and the index are independent of the choice of principal homogeneous space structure on $V,$ so that our notation is
not really so bad.}  One shows easily that for all $V$, the period divides the index and the two quantities have
the same prime divisors.  (These and other elementary facts about the period and the index are reviewed in Section 
2.)  To say more is the \textbf{period-index problem} for the Weil-Ch\^atelet
group $H^1(k,A)$. 
\\ \\
The classic paper on the period-index problem in the context of abelian varieties is
[Lang-Tate].  Almost all subsequent work ([Lichtenbaum], [Cassels], [Stein], [O'Neil]) has focused on the 
case where $A = E$ is an elliptic curve.  In this case, if one restricts to classes in
$H^1(k,E)[p]$ -- i.e., to genus one curves $C/k$ with prime period $p$ -- there is a simple dichotomy: either
$C$ acquires a rational point over a degree $p$ field extension, or it does not, in which case it necessarily
acquires a rational point over a field extension of degree $p^2$.  Deciding which occurs is a basic problem in
the arithmetic of genus one curves.  Some important results in this direction are as follows.
\\ \\
([Lichtenbaum], [O'Neil]) If $k$ has vanishing Brauer group, then period equals index for all genus one curves over $k$.  
\\ \\
([Lichtenbaum]) The period equals the index for all genus one curves over a $p$-adic field.
\\ \\
([Cassels], [O'Neil]) If $C/k$ is a genus one curve over a number field which has rational points everywhere locally, then
the period equals the index.
\\ \\
Cassels also gave examples of genus one curves over $\Q$ of period $2$ and index $4$.  A generalization was proved in [Clark]:
let $k$ be a number field, $p$ a prime number and $E/k$ an elliptic curve with full $p$-torsion defined over $k$.  Then there
exists an infinite subgroup of $H^1(k,E)[p]$ all of whose nonzero elements have index $p^2$.
\\ \\
Rather fewer results are known for abelian varieties.  Most relevant to our purposes are the following:
\\ \\
([Lang-Tate]) Let $A/k$ be an abelian variety over a $p$-adic field.  Suppose that $A$ has good reduction, and let $n$ be prime
to $p$.  Then any element $V$ of $H^1(k,A)[n]$ has index $n$: indeed, the splitting fields $l$ of $V$ are precisely those 
for which the ramification index $e(l/k)$ is divisible by $n$.
\\ \\
In the case where $k$ is a $p$-adic field and $A/k$ is an abelian variety admitting an analytic uniformization, there is a 
similarly
complete account of splitting fields of principal homogeneous spaces due to [Gerritzen].  Later we shall give precise
statments, and even proofs, of some of his results.
\\ \\
The goal of the present work is twofold.  First, we examine the period-index problem in the higher-dimensional case and derive
analogues of the aforementioned results for elliptic curves, including the result of [Clark].  This latter result is
then applied to the horizontal variation of the Shafarevich-Tate group of a principally polarized abelian variety.  Here are the
statements of the main theorems.
\begin{thm}
Let $g$ be a positive integer and $p$ be a prime number.  Let $k$ be any one of the following fields: \\
$\bullet$ a sufficiently large $p$-adic field $k = k(p,g)$; \\
$\bullet$ a sufficiently large number field $k = k(g)$; or \\
$\bullet$ the maximal unramified extension $k_{\infty}$ of a $p$-adic field $k$ containing $\Q(\mu_p)$. \\
Then there exists an abelian variety $A/k$ and a principal homogeneous space $V \in H^1(k,A)[p]$
of index at least $p^g$.
\end{thm}
\noindent
Remark: 1) The $p$-adic case is essentially a theorem of Gerritzen, and the remaining cases are deduced from it.
I was not aware of
Gerritzen's work until quite recently -- indeed, an earlier draft of this paper contained an independent 
proof of a weaker result, namely with the index replaced by the \emph{Galois} index (see Section 2.3), from which we deduced
(only) that the minimal degree of a splitting field of a principal homogeneous space of prime period could be arbitrarily large.
\\ \indent 
2) In the $p$-adic case, sufficiently large means that $k$ contains the $p$th roots of unity and
$[k:\Q_p] \geq g-2$.  In the number field case it is a bit more complicated, and we address the question of
making $k$ as small as possible only for $p = 2$, $g \leq 3$ (in which case we show that we may take $k = \Q$). \\
\indent
3) The maximal unramified extension of $\Q_p$ has trivial Brauer group; hence the theorem of [Lichtenbaum] and
[O'Neil] does not generalize in the form ``period equals index'' for higher-dimensional abelian varieties (but see Theorem 3). \\
\indent
4) As we shall see, an easy restriction of scalars argument allows one to deduce that over
any local or global field, there exist principal homogeneous spaces of prime period and arbitarily
large index.
\begin{thm}
Let $A/k$ be a $g$-dimensional principally polarized abelian variety over a 
$p$-adic field, and let $V \in H^1(k,A)[n]$.  Suppose that \textbf{either} \\
$\bullet$ $n$ is odd;  \textbf{or} \\
$\bullet$ $A[n]$ is a Lagrangian $\gk$-module. \\
Then $V$ is split by a field extension of degree at most $(g!)n^g$.
\end{thm}
\noindent
Remark: The definition of a Lagrangian Galois-module structure on $A[n]$ is 
given in Section 6.3, but note that at least the trivial Galois module structure is Lagrangian.  
\begin{thm}
Let $A/k$ be a $g$-dimensional strongly principally polarized abelian variety and $V \in H^1(k,A)[n]$ a principal
homogeneous space.  Suppose that \emph{either} \\
$\bullet$ $k$ has trivial Brauer group; \textbf{or} \\
$\bullet$ $k$ is a number field and $V \in \Sha(A,k)[n]$ is a locally trivial class.  \\
Then $V$ is split over a field
extension of degree at most $(g!)n^g$.
\end{thm}
\noindent
There is an immediate application to an upper bound for the \textbf{visibility dimension}
of a class in $\Sha(A/k)[n]$.  We refer the reader to [Agashe-Stein] and [Cremona-Mazur] for a
discussion of visibility of principal homogeneous spaces.  We recall only the following fact:
if $V \in H^1(k,A)$ can be split by a field extension of degree at most $N$, then 
the visibility dimension of $V$ is at most $g \cdot N$, where $g = \dim V$ [Agashe-Stein, Prop. 1.3].
Thus Theorem 3 gives the following upper bound for the visibility dimension
of a locally trivial class (compare with [Agashe-Stein, Prop. 2.3 and Remark 2.5]). 
\begin{cor}
Let $\eta \in \Sha(A/k)[n]$, where $A/k$ is a $g$-dimensional strongly principally
polarized abelian variety over a number field $k$.  Then the visibility dimension of
$\eta$ is at most $g \cdot (g!) \cdot n^g$.
\end{cor}
\noindent
Remark: Because there could be many other ways to visualize $V$, there is no reason to
believe that the upper bound of Corollary 4 is sharp.  Indeed, [Mazur] showed that
any nontrivial element of $\Sha(E/k)[3]$ has visibility dimension $2$, and [Klenke] showed
the same result for all elements of $H^1(k,E)[2]$ (i.e., even those with index $4$).  
\\ \\
So far we have been pursuing results in analogy with the elliptic curve case.  However, an interesting
new phenomenon emerges starting in dimension two: the period-index problem in the Brauer group
of $k$ begins to an exert an influence on the period-index problem for WC-groups.
\begin{thm}
Let $n$ and $a$ be positive integers and $k$ a field such that every class in
$Br(k)[n]$ can bes split over a field extension of degree dividing $n^a$.  (In particular,
if $k$ is a local or global field we may take $a = 1$ for all $n$.)  Let $A$ be a strongly
principally polarized abelian variety, and assume \textbf{either} of the following:  \\
$\bullet$ $\car(k) \neq 2$ and $n$ is odd; \textbf{or} \\
$\bullet$ $A[n]$ is a Lagrangian $\gk$-module. \\
Then any $V \in H^1(k,A)[n]$ can be split over a field extension of degree at most $(g!)n^{a+g}$.
\end{thm}
\noindent
Remark: The optimal bound for the degree of a splitting field of
$V \in H^1(k,A)[n]$ is $n^{2g}$.  Lang and Tate showed that this bound could be attained
over iterated Laurent series fields (the precise result is recalled in Section 2).  Theorem 5
shows that it is indeed necessary to consider such complicated fields in order to get large odd indices.
\begin{thm}
Let $p$ be a prime number and a $A/k$ be a strongly principally polarized abelian variety over a 
number field.  Assume that both $A[p]$ and $NS(A)$ (the N\'eron-Severi group of $A$) are
trivial as Galois modules.  Then there exists an infinite subgroup $G \subset H^1(k,A)[p]$
such that every nonzero $V \in G$ has index exceeding $p$.
\end{thm}
\begin{thm}
Let $A, \ ,k, \ p$ be as in the statement of the previous theorem.  Then for every positive
integer $a$ there exist infinitely many degree $p$ field extensions $l/k$ such that
$\# \Sha(A/l)[p] \geq p^a$.
\end{thm}
\begin{cor}(Horizontal variation of $\Sha$)
Let $A/k$ be any $g$-dimensional principally polarized abelian variety over a number field and
$p$ any prime number.  There exists a function $F(g,p)$ such that
\[\sup_{l/k \ : \ [l:k] \leq F(g,p)} \#(A/l)[p] = \infty. \]
For instance, one can take $F(g,p) = p \cdot 2^{2g} \cdot \# GSp_{2g}(\F_p) \cdot \#GL_{4g^2}(\F_3)$.
\end{cor}
\noindent
Remark: It seems likely that one can take $F(g,p) = p$.
\\ \\
The organization of the paper is as follows.  Section 2 collects some preliminary results on
period-index problems in general, and in the Weil-Ch\^atelet group in particular.
\\ \\
In Section 3 we recall Gerritzen's work on analytically uniformized abelian varieties; this work
is then used to give a proof of Theorem 1.  An additional argument using modular curves
is used to get $3$-dimensional principal homogeneous spaces over $\Q$ of period $2$ and index $8$.
\\ \\
In Section 4 we formulate two separate period-index problems for any variety $V/k$, one defined
in terms of the Albanese variety and the other in terms of the Picard variety.  If $V$ is a curve,
then of course these are one and the same problem; in the higher-dimensional case, this is a convenient
framework for transferring results about divisors to results about zero-cycles.
\\ \\
The core of the paper is Section 5, in which we consider an adaptation of O'Neil's \textbf{period-index obstruction}
to our higher-dimensional context.  The basic theta group construction can be done verbatim
using abelian varieties instead of elliptic curves, but complications arise in the higher-dimensional
case.  Indeed, since the index of a variety involves its zero-dimensional geometry (least
degree of an effective $k$-rational zero cycle) and the period-index obstruction involves
its codimension-one geometry (the obstruction to a rationally defined divisor class coming
from a rational divisor), when $\dim V > 1$ seems \emph{a priori} surprising that the
non/vanishing of the obstruction map $\Delta$ should be related to the index.  Nevertheless,
using the Albanese/Picard formalism, we show that as in the one-dimensional case, consequences can be drawn both from the 
vanishing and the non-vanishing of $\Delta$ (at least under certain geometric hypotheses on $A/k$ that are always satisfied
in the one-dimensional case).  We deduce Theorem 3 immediately from this setup, 
as in [O'Neil].  
\\ \\
In [Clark] we used Mumford's theory of the Heisenberg group to get an ``explicit'' form
of the period-index obstruction map when $A = E$ has full level $n$-structure.  In Section 6 we
look at the Galois cohomology of Heisenberg groups in more detail.  The results
directly imply Theorem 5 and are used in the proofs of the remaining theorems.
\\ \\
In Section 7 we give the proof of Theorem 6, and in Section 8 we prove Theorem 7 and Corollary 8, following 
roughly the same strategy as in the one-dimensional case, which was treated in [Clark].
\\ \\
We end in Section 9 with some comments.
\\ \\
Acknowledgements: This work would not be possible without the basic insights provided
by William Stein and Catherine O'Neil.  Many others have made valuable contributions
along the way: it is a pleasure to thank Kevin Buzzard, Romyar Sharifi, Florian Herzig, 
Joost van Hamel, E. Z. Goren and Martin Hillel Weissman.
\\ \\
Notation and conventions:
Throughout this paper $k$ denotes a field, assumed (mostly for convenience) to be perfect, 
$\overline{k}$ denotes some fixed algebraic closure of $k$, and $\gk = \Gal(\overline{k}/k)$
is the absolute Galois group of $k$.  When we say
$V/k$ is a variety, we mean that it is a smooth, projective and geometrically irreducible
$k$-variety.
\\ \\
The letters $n$ and $p$ shall always denote respectively a positive integer and a prime
number, and they are always assumed to be \emph{indivisible} by the characteristic of $k$.
Related to this, we let $\mu_n$ and $A[n]$ denote the $n$-torsion subgroups of
$\Gm$ and $A/k$ an abelian variety respectively.  These are finite flat $k$-group
schemes which, owing to our assumption on $n$, are \'etale.  They may thus be viewed
as $\gk$-modules.  
\\ \\
If $M$ is a $\gk$-module, we write $H^i(k,M)$ for the Galois cohomology group
$H^i(\gk,M)$.  If $l/k$ is a field extension, we write $H^i(l/k,M)$ for the kernel
of the natural (restriction) map $H^i(k,M) \ra H^i(l,M)$.  
\\ \\
If $M$ is a finite $\gk$-module, we write $M^* = \Hom(M,\Q/\Z(1)) = \Hom(M,\Gm[\tors])$
for its Cartier dual and $M^{\vee} = \Hom(M,\Q/\Z)$ for its Pontrjagin dual.
\\ \\
If $M$ is a trivial $\gk$-module, we often speak of elements of $H^1(k,M) =
\Hom(\gk,M)$ as ``characters,'' even though the extension cut out
by such a class need not be cyclic.  The point is that such classes behave in a very
easily understood way: they have a \emph{unique} minimal splitting field, which is an
abelian extension of $k$.
\section{Period-index problems in Galois cohomology}
\noindent
The first two sections collect foundational material on period-index problems in
general and in Weil-Ch\^atelet groups in particular.  None of these results are due
to the present author.  The third section contains some discussion of possible
variants on the index; except for the definition of the ``mindex,'' it may safely
be omitted on a first reading.
\subsection{The period and index of a cohomology class}
We insist on presenting the first few results in ``extreme generality'': namely,
let $M$ be any commutative $\gk$-module, $i > 0$ an integer and $\eta \in H^i(k,M) =
H^i(\gk,M)$ a Galois cohomology class.  We will define the period and index for
$\eta$.  On the one hand, we do this because the period-index problem \emph{is}
interesting for Galois modules other than $A(\overline{k})$: e.g., the case of
$H^i(k,G)$ for $G/k$ a commutative algebraic group includes both the present
case and the case of the Brauer group, $Br(k) = H^2(k,\Gm)$, where indeed much
more work has been done than in WC-groups.\footnote{In an earlier draft, the Brauer group
received equal billing with WC-groups.  Since a large part of what we had to say about the
Brauer group was expository in nature, we have reworked things to force $Br(k)$ to play
a subsidiary role.}  On the other hand, by working in such generality we show that there is
``nothing up our sleeves'': if one works in the context of $Br(k)$ or $H^1(k,A)$ then one
might try to use aspects of the theory of division algebras or abelian varieties in order
to prove the results, and this would be working too hard.
\\ \\
For $i > 0$, we define two
numerical invariants of a class $\eta \in H^i(k,M)$.  Recall that $H^i(k,M)$ is an abelian
torsion group.  Thus we may define the \textbf{period}
$n(\eta)$ to be the order of $\eta$ as an element of $H^i(k,G)$.   The \textbf{index}
$i(\eta)$ is the greatest common divisor of all degrees of finite field extensions $l/k$ that split
$\eta$, i.e., such that $\eta|_{\mathfrak{g}_L} = 0$.  
\begin{prop}
Let $\eta \in H^i(k,M)$ be any Galois cohomology class, with $i > 0$. \\
a) The period $n(\eta)$ divides the index $i(\eta)$. \\
b) The period and index of $\eta$ have the same prime divisors. \\
c) (Non-reduction) If $l/k$ is a field extension of degree prime to $n(\eta)$,
then $n(\eta|_l) = n(\eta)$ and $i(\eta|_l) = i(\eta)$. \\
d) (Primary decomposition) Let $\eta = \eta_1 + \ldots + \eta_r$
be the primary decomposition of $\eta$ corresponding to the
factorization $n(\eta) = p_1^{a_1} \cdots p_r^{a_r}$, i.e.,
$\eta_i = \frac{n}{p_i^{a_i}} \eta$.  Then $i(\eta) = \prod_{i=1}^r i(\eta_i)$.
\end{prop}
\noindent
Proof: For part a): if $l/k$ is a degree $n$ splitting field for $\eta$, then
$0 = Cor(Res(\eta)) = n \eta$ (see [CL, Prop. 7.6]), showing that the period 
divides the index. \\ \indent
For part b): let $p$ be a prime number which is prime to the period of $\eta$.  Let
$f/k$ be a finite Galois splitting field, and let $l/k$ be the subextension corresponding to some
Sylow $p$-subgroup of $\ggg_{f/k}$.  Consider the restriction of $\eta$ to $\ggg_l$.  On the one hand,
its period divides the period of $\eta$ hence remains prime to $p$, but on the other hand
by construction this class is split over a field extension whose degree is a power of $p$, so by the
first part of the proposition its period is a power of $p$.  Therefore $l$ is a splitting field
for $\eta$, so the index of $\eta$ is prime to $p$.  \\ \indent The remaining parts are 
routine.
\begin{prop}
Let $M$ be a finite $\gk$-module and $\eta \in H^1(k,M)$. \\
a) The class $\eta$ can be split by a field extension of degree at most $\#M$, so
$i(\eta) | \# M$. \\
b) If $M$ is a trivial $\gk$-module, the index of $\eta$ is attained: there
exists $l/k$ of degree $i(\eta)$ such that $\eta_|{\ggg_l} = 0$.
\end{prop}
\noindent
Proof (Lenstra): Let $\xi: \gk \ra M$ be a one-cocycle representing $\eta$.  Define
$H \subset \gk$ to be the subset of elements $\sigma$ such that $\eta(\sigma) = 0$.
Despite the fact that $\xi$ is not necessarily a homomorphism, one nevertheless has
that $H$ is a subgroup of $\gk$ and that $\xi$ induces an injective map of sets
\[\xi: \gk/H \hookrightarrow M, \]
where $\gk/H$ is the right coset space.  Thus $H = \ggg_l$ corresponds to a
splitting field of degree at most $\#M$, giving the first statement of part a).  In particular 
$i(\eta) \leq \#M$.  If $M$ has prime-power order, then we must have
that  $i(\eta) | \#M$, and the general case follows by primary decomposition.  For part b),
if $M$ is trivial, then $H^1(k,M) = \Hom(\gk,M)$ is just the group of ``$M$-valued characters
of $\gk$.''  In particular, every class $\eta \in H^1(k,M)$ has a unique minimal splitting field,
namely the fixed field $l$ of  $\ker(\eta)$, an abelian extension.
In this case $\ggg_{l/k} \hookrightarrow M$ is a homomorphism of groups and part b) follows.
\\ \\
Remark: The hypothesis $i = 1$ in the proposition is necessary: e.g., elements of
$Br(k)[n] = H^2(k,\mu_n)$ need not have index $n$.
\subsection{Results on Weil-Ch\^atelet groups}
\begin{cor}
Let $A/k$ be an abelian variety of dimension $g$, $n$ a positive integer -- indivisible, as always,
by the characteristic of $k$ -- and $V \in H^1(k,A)[n]$.  Then $V$ is split by some field
extension $l/k$ of degree at most $n^{2g}$; in particular $i(V)$ divides $n^{2g}$.  If $A[n]$
is $\gk$-trivial, $V$ is split by a field extension of degree dividing $n^{2g}$.
\end{cor}
\noindent
Proof: For any field extension
$l/k$ we have a Kummer sequence
\[0 \ra A(l)/nA(l) \ra H^1(l,A[n]) \ra H^1(l,A)[n] \ra 0, \]
compatible with the restriction maps from $k$ to $l$.  Thus it is sufficient to trivialize
any lift $\xi$ of $V$ to $H^1(k,A[n])$.  Since $\#A[n](\overline{k}) = n^{2g}$, the conclusion
follows from Proposition 10.
\\ \\
If there is no restriction on the field $k$, the bound $i \ | \ n^{2g}$ is optimal:
\begin{prop}([Lang-Tate])\\
a)
Suppose we have a field $k$, a $g$-dimensional abelian variety $A/k$ and a positive integer $n$ satisfying 
the following hypotheses:
\begin{itemize}
\item{$A[n]$ is $\gk$-trivial.}
\item{For every finite extension $l/k$, $A(l)$ is $n$-divisible.}
\item{There exists a Galois extension $l_0/k$ with Galois group isomorphic to $(\Z/n\Z)^{2g}$.}
\end{itemize}
Then for all $a$, $1 \leq a \leq 2g$, there is a cohomology class $\eta \in H^1(k,A)$ of period $n$ and index $n^{a}$. 
\\
b) The hypotheses of part a) are satisfied for an isotrivial abelian variety (i.e., one arising
by basechange from $\C$) over the iterated
Laurent series field $k_{2g} := \C((t_1))\cdots ((t_{2g}))$.
\end{prop}
\noindent
We will also need the following simple fact:
\begin{prop}([Lang-Tate])
Let $v$ be a discrete valuation on a field $k$, with completion $k_v$.  Let $n$ be a positive integer and
$A/k$ be an abelian variety with $A[n](k) = A[n](\overline{k})$.  Suppose moreover that $k$ contains
the $n$th roots of unity.  Then the local 
restriction map $H^1(k,A)[n] \ra H^1(k_v,A)[n]$ is surjective.
\end{prop}
\noindent

\subsection{Variations on the index}
If $k$ were not perfect, the above definition of the index would allow for inseparable field extensions.  
But then it is natural to wonder what happens when one
restricts to separable splitting fields, leading to the notion 
of the \textbf{separable index} $i_s(\eta)$.  The separable index has all
the formal properties of the index: indeed, if $i$ is replaced by $i_s$ and $\overline{k}$ by 
$k^{\sep}$, all the results of the paper
remain valid.  Moreover, it seems likely that $i(V) = i_s(V)$
for all principal homogeneous spaces $V$ -- in the terminology of
Section 4, this is true for the Picard index of any variety and in particular in dimension one 
[Lichtenbaum], [Harase].  My feeling that these
issues of separability are of relatively minor interest is the source of the 
assumption on the perfection of $k$.  The reader who feels otherwise
is invited to consider the general case.
\\ \\
The Galois index: In the definition of the index does of course allow not-necessarily
normal field extensions $l/k$.  The \textbf{Galois index} $i_G(\eta)$, defined using this
restriction, can indeed properly exceed the index.  The question of equality
$i_G(\eta) = i(\eta)$ for classes in the Brauer group is equivalent to the question of
which division algebras are crossed-product algebras; by deep work of [Amitsur]
we know that the answer to this question is -- for most periods $n$ --  generically negative.    \\ \\
Here is an example to show that $i_G > i$ can occur in WC-groups.  Let $k$ be a $p$-adic
field with residue cardinality $p^a$ and $A/k$ an abelian variety with good reduction.  Suppose that
$\ell > p^a$ is a prime and $V \in H^1(k,A)[\ell]$ is a nontrivial element.  Then $i_G(V) > i(V)$.
This is obtained by repeated application of an aforementioned result of Lang and Tate, which can be
rephrased as: let $m/f$ be an extension of $p$-adic fields and $A/f$ an abelian variety with
good reduction.  Then the natural map
$H^1(f,A)[\ell] \ra H^1(m,A)[\ell]$ is the zero map if $\ell \ | \ e(m/f)$ and is injective otherwise.
In particular, $i(V) = \ell$.  But let $l/k$ be any Galois splitting field for $V$.  We claim that $l$ must contain $\mu_{\ell}$, 
hence $\ell - 1 = \ [k(\mu_{\ell}):k] \ | \ [l:k]$, so that
$\ell - 1$ divides $i_G(V)$.  Indeed, $l/k$ can be decomposed as a tower $k \subset l_1 \subset l_2 \subset l_3 = l$,
where $l_1/k$ is unramified, $l_2/l_1$ is totally ramified of degree prime to $p$, and $l_3/l_2$ is totally ramified and of degree 
a power
of $p$.  We get that the restriction maps from $l$ to $l_1$ and from $l_2$ to $l_3$ are both injective on elements of
period $\ell$, so $V|_{l_1} \neq 0$ and $V|_{l_2} = 0$.  That is, $V|_{l_1}$ is killed by the totally tamely ramified
extension $l_2/l_1$, which is necessarily of the form $l_2 = l_1[T]/(T^a-\pi)$, where $\ell \ | a$
and $\pi$ is a uniformizer of $l_1$.  But if this extension is Galois, $l_1$ contains the $\ell$th roots of unity,
establishing the claim.  Of course we assumed that such a nontrivial
$V$ exists, which is not always the case (indeed, for fixed $A/k$ such a $V$ exists for at most
finitely many primes $\ell$), but we may certainly arrange for such classes to exist: e.g.,
fix a prime $p$, and let $\ell$ be a prime such that
$p+ 1 < \ell < p+1+\sqrt{2}p$.  By the Hasse-Weil-Waterhouse theorem, there is an (ordinary) elliptic 
curve $E/\F_p$ such that $E(\F_p) \cong \Z/\ell \Z$; let $\tilde{E}/\Q_p$ be its canonical lift.
Then $\Z/\ell\Z \cong \tilde{E}(\Q_p)/\ell \tilde{E}(\Q_p) \cong H^1(\Q_p,E)[\ell]$.  The last
isomorphism is by Tate's duality theorem.  
\\ \\
This example should be compared with the case of principal homogeneous spaces $V$ over $p$-adically uniformized
abelian varieties discussed in the next section: for all such $V$ we have $i(\eta) = i_G(\eta)$.  It would be nice
to know whether $i = i_G$ in the case of good reduction and period equal to the residue characteristic.
\\ \\
Finally and most importantly, one can ask for examples in which the index is not attained, i.e., such that the greatest common 
divisor of
degrees of splitting fields is not itself the degree of a splitting field.  I do not know of such an example (anywhere in 
Galois cohomology). 
It is well-known that the index is attained for classes in the Brauer group $H^2(k,\Gm)$.  The attainment of
the index for elements of the Weil-Ch\^atelet group of an elliptic curve was observed in [Lang-Tate].  In contrast, the 
attainment of the index is an important \emph{open problem} in the Weil-Ch\^atelet group of a higher-dimensional abelian 
variety, and our ignorance of this attainment leads us to define the \textbf{mindex} $m(\eta)$ of a Galois cohomology class
as the minimal degree of a splitting field.  
\\ \\
Clearly the best response to the possibility of $i(V) < m(V)$ is to give upper bounds on the mindex
and lower bounds on the index.   Fortunately enough, it turns out that most of our main results
results are phrased in this way.
\section{Large indices over local, strictly local and global fields}
\subsection{Travaux de Gerritzen}
In this section we give an account of some work of Gerritzen on the period-index problem 
in the Weil-Ch\^atelet group of an analytically uniformized abelian variety.  Because these results are closely
related to the proof of Theorem 1, we will give complete proofs.
\begin{prop}(Gerritzen)
Let $k$ be any field, with absolute Galois group $\gk$.  Let $\tilde{A}$ be a $\gk$-module and
$\Gamma \subseteq \tilde{A}$ a $\gk$-submodule which is torsionfree as a $\Z$-module and such
that $\Gamma^{\gk} = \Gamma$; put $A := \tilde{A}/\Gamma$.  
Suppose also that $H^1(l,\tilde{A}) = 0$
for all finite extensions $l/k$.  Let $\eta \in H^1(k,A)$ be a class of
exact period $n$.  Then: \\
a) $\eta$ has a \emph{unique} minimal splitting field $L = L(\eta)$. \\
b) The extension $l/k$ is abelian of exponent $n$. \\
c) $i(\eta) \ | \ n^g$, where $g = \dim_{\Q}(\Gamma \otimes \Q)$ is the rank of
$\Gamma$.
\end{prop}
\noindent
\newcommand{\gl}{\mathfrak{g}_l}
Proof: We take $\gk$-cohomology of the short exact sequence
\begin{equation}
0 \ra \Gamma \ra \tilde{A} \ra A \ra 0,
\end{equation}
and using $H^1(k,\tilde{A}) = 0$, we get an injection
$\delta: H^1(k,A) \hookrightarrow H^2(k,\Gamma)$.  Also,
since $\overline{\Gamma} = \Gamma \otimes \Q$ is cohomologically
trivial, there is a canonical isomorphism
$H^2(k,\Gamma) \cong H^1(k,(\Q/\Z)^g) = H^1(k,\Q/\Z)^g =: X_A(\gk)$,
i.e., $g$ copies of the character group of $\gk$.  Indeed,
the assumptions are such that for any finite extension $l/k$,
we may view $(1)$ as a sequence of $\mathfrak{g}_l$ modules
and get the same result: we get for every $l$ an injection
$\delta_l: H^1(l,A) \hookrightarrow H^1(l,\Q/\Z)^g = X_A(\gl)$,
and these various maps are compatible with restriction.  Thus splitting
$\eta$ is equivalent to splitting the character $\delta_k(\eta)$.
But $\delta_k(\eta)$ cuts out a field extension $l(\eta)/k$ which
is abelian of exponent $n$, of order dividing $n^g$, and evidently the
unique minimal splitting field, completing the proof.
\\ \\
Recall that a $g$-dimensional abelian variety defined over a $p$-adic
field admits an analytic uniformization if it is isomorphic, as a
rigid $k$-analytic group, to $\Gm^g/\Gamma$, where
$\Gamma \cong \Z^g$ is a discrete subgroup.
\begin{thm}(Gerritzen)
Let $A/k$ be a $g$-dimensional analytically uniformized abelian variety over a $p$-adic field. \\
a) For any $V \in H^1(k,A)[n]$, $i(V) \ | \ n^g$. \\
b) If $\Gamma \subseteq nA(k)$, then
$\delta_k: H^1(k,A)[n] \ra H^1(k,\Z/n\Z)^g$ is an isomorphism of finite groups.
\end{thm}
\noindent
Proof: Part a) follows immediately from the Proposition, by taking
$\tilde{A} = \Gm^g(\overline{k})$: note that $H^1(l,\tilde{A}) = 0$
for all $l/k$ by Hilbert 90.  As for part b), we have an injection
of finite groups, so it's enough to see that they have the same
cardinality.  We recall the following important duality theorem of 
Tate: for all finite $\gk$-modules $M$, $H^1(k,M^*) \cong
H^1(k,M)^{\vee}$ [CG, $\S$ II.5.2, Theorem 2].  Applying this to $M = \Z/n\Z$, we get that
$H^1(k,\Z/n\Z) \cong H^1(k,\mu_n) \cong \# k^*/k^{*n} \cong
(\Z/n\Z)^r$ for some positive integer $r$.  Thus
\[\#H^1(k,\Z/n\Z)^g = n^{rg}. \]
On the other hand, by Tate's local duality theorem for abelian varieties [Tate],
\[H^1(k,A)[n] \cong \frac{A(k)}{nA(k)} =
\frac{\Gm^g(k)/\Gamma}{(\Gm^g(k))^n/(\Gamma \cap
\Gm^g(k)^n)} \cong \left(\frac{k^*}{k^{*n}}\right)^g, \]
the last isomorphism because of the assumption that every
element of $\Gamma$ is an $n$th power in $\Gm^g(k)$.  (Note that
without this
assumption on $\Gamma$, the map $\delta_k$ clearly does not have
to be an isomorphism.)
Thus
we also have $\#H^1(k,A)[n] = n^{rg}$, completing the proof.
\subsection{The proof of Theorem 1 in the $p$-adic case}
We now specialize to the following situation: let 
$k/\Q_p$ be a finite extension of degree $a$; we assume that
$k$ contains the $p$th roots of unity.  One knows that
\begin{equation}  
[k^*:k^{*p}] = p^{a+2}.
\end{equation}
Let $E/\Q_p$
be an elliptic curve with analytic uniformization 
$\Gm / \langle q \rangle$, where $q = p^p$ (say) is
a $p$th power in $\Q_p$, \emph{a fortiori} in $k$ -- this choice
is so that the hypothesis of Theorem 15b) is satisfied for $n =p$.
As we saw in the proof of Theorem 15, $\dim_{\F_p} H^1(k,E)[p] = 
p^{a+2}$.  Let $P_1, \ldots, P_{a+2}$ be an $\F_p$-basis
for $H^1(k,E)[p]$.  Then Proposition 14 associates to each
$E_i$ a unique minimal splitting field $l_i$, such
that $l_i/k$ is cyclic of degree $p$.  Because of the
injectivity of $\delta_k$, the characters $\{\delta_k(P_i)\}_{i=1}^{a+2}$,
remain $\F_p$-linearly independent, so cut out an abelian extension
$l^{(p)} = l_1 \cdots l_{a+2}$ of exponent $p$.  
Let $P = P_1 \times \cdots \times P_{a+2}$, viewed as a principal homogeneous
space over the (analytically uniformized) abelian variety $A = E^{a+2}$.
Evidently a field $l$ splits $P$ exactly when it splits each
$P_i$, so $l^{(p)}$ is the unique minimal splitting field of $P$,
which therefore has index $p^{a+2} = p^{\dim A}$.
\subsection{Applications to $k_{\infty}$}
Let $k/\Q_p(\mu_p)$ be a $p$-adic field containing the $p$th roots of unity, and
denote by $k_{\infty}$ the maximal unramified extension of $k$.  Let $l_1/k$ be the
unique unramified extension of degree $p$, and let $l_1, l_2, \ldots, l_{a+2}$
be an $\F_p$-basis for the set of abelian $p$-extensions of $k$ (i.e., a linearly disjoint
set of extension fields whose compositum is the maximal abelian extension of exponent
$p$).  Under the bijection $\delta_k$ the $l_i$'s are splitting fields
of homogeneous spaces $P_1, \ldots, P_{a+2}$; let $Q := P_2 \times P_3 \times \cdots \times P_{a+2}$.
\\ \\
We claim that the index of $Q/k_{\infty}$ is the same as its index over $k$, namely $p^{a+1}$.  Indeed, let
$m/k_{\infty}$ be a degree $N$ splitting field of $Q/k_{\infty}$.  Then $Q/k$ is split by some finite extension
$m'/k$, such that $m' \cap k_{\infty} = k_b$ (the unramified extension of degree $b$)
for some $b$, and $[m':k_b] \leq N$.  By construction, the unique minimal splitting field for $Q$ is disjoint
from $k_b$, so the index of $Q$ is not reduced by restriction to $k_b$.  Thus
$p^{a+1} \ | N$, and we conclude that the index of $Q/k_{\infty}$ is $p^{a+1}$.  Since we may view
$k_{\infty}$ as the maximal unramified extension not just of $k$ but of any $k_b$, we may arrange for
$a+1 = [l:\Q_p] \geq g$, completing the proof.
\subsection{Applications to number fields}
Fix $p$ and $g$.  Let $E/\Q$ be any elliptic curve with potentially multiplicative reduction at $p$ (e.g.,
one with $j$-invariant $\frac{1}{p}$).
Choose a field extension $k/\Q$ which is sufficiently large in the following senses: \\
a) $E[p]$ is a trivial $\gk$-module ($\implies \Q(\mu_p) \subset k$), \\
b) $k$ has a place $v|p$ of local degree at least $g-2$ \\
c) $E/k_v$ has split multiplicative reduction. \\ \\  
From Section 3.2, we have a class $\eta_v \in H^1(k_v,E^g)$ with period $p$ and index $p^g$;
by Proposition 13 there exists a class $\eta \in H^1(k,E^g)[p]$
mapping to $\eta_v$.  Thus $\eta$ is a global class of period $p$ and index \emph{at least} $p^g$.
\\ \\
\newcommand{\Res}{\operatorname{Res}}
Examples over $\Q$: Let $A/k$ be as above and consider $\Res_{k/\Q} A$, the abelian variety
over $\Q$ obtained from $A/k$ by restriction of scalars (or ``Weil restriction'').  There is a canonical
isomorphism $H^1(k,A) = H^1(\Q,\Res_{k/\Q} A)$ (cf. [Agashe-Stein]); if $\eta \in H^1(k,A)[p]$ has index
$p^i$, then the corresponding class $\eta^{\Q}$ in $H^1(\Q,\Res_{k/\Q} A )$ must have period $p$ and index
at least $p^i$: indeed, if $m/\Q$ is a splitting field for $\eta^{\Q}$, then $ml/l$ splits $\eta$.  Notice however
that $\dim \Res_{k/\Q} A = [k:\Q] \dim A$, so in general we do not get classes of period $p$ and index
at least $p^{\dim A}$ over $\Q$ using this trick.
\\ \\
But let us look at the case $p = 2, \ g = 3$.\footnote{\emph{A fortiori} the construction works for
$g = 2$.  Indeed, when $g = 2$ one may take the period to be prime to $p$; cf. [Lichtenbaum].}
\begin{prop}
There exists an abelian 3-fold $A/\Q$ and a principal homogeneous space $X \in H^1(\Q,A)$
of period $2$ and index at least $8$.
\end{prop}
\noindent
\emph{Proof}: From the above discussion, we need only find an elliptic curve $E/\Q$ with
$E[2] = E[2](\Q)$ and with split multiplicative reduction at $2$.  To do this:
\\ \\
On the one hand, the elliptic curve labeled 4290Z2 in Cremona's 
tables,
with minimal Weierstrass model
$E: y^2 + xy = x^3 - 66x$ fits the bill.\footnote{The time when Kevin Buzzard came up with this equation before my eyes
is a five-minute interval that I shall not soon forget.}
Alternately, we have the
following
\begin{prop}
Let $X = X(\Gamma)$ be the modular curve corresponding to the congruence subgroup 
$\Gamma := \Gamma(2) \cap \Gamma_1(4)$ of $SL_2(\Z)$.  Then there are infinitely many
$P \in X(\Q)$ with split multiplicative reduction at $2$.
\end{prop}
\noindent
\emph{Proof}: One knows that $X/\Q$ is the fine moduli space for $\Gamma$-structured elliptic curves over
$\Q$, and also that it has genus $0$.  As a smooth genus $0$ curve, $X$ is isomorphic to 
$\mathbb{P}^1$ over $\Q$ if and only if it has points everywhere locally, or even at every finite
place $p$.  But the Tate curve $\Gm/<p^4>$ gives a point on $X(\Q_p)$.  So $X \cong_{\Q} \mathbb{P}^1$ (as is well-known).
Now let $P_0 \in X(\Q_2)$ be a ``Tate point'' as above.  We claim that there is an analytically
open neighborhood $U$ of $P_0$ in $X(\Q_2)$ consisting of points with Tate uniformizations -- i.e. with
split multiplicative reduction.  Indeed, consider the universal $\Gamma$-structured elliptic curve
$\mathcal{E}$ over an analytic neighborhood of $P_0$ -- it is given by a Weierstrass equation with
analytically varying coefficients.  The condition that the $\Gamma$-structured elliptic curve
over $P$ have split multiplicative reduction is that the quantity $c_4$ be nonzero 
and that the quadratic (lowest-degree) form of the Weierstrass equation factor over
$\Q_2$ -- these are both open conditions.  Since $X(\Q)$ is dense in $X(\Q_2)$, $X(\Q) \cap U$
is infinite, which was to be shown.
\section{Albanese and Picard varieties}
\subsection{Background on Albanese, Picard and N\'eron-Severi}
This section contains foundational material on Albanese varieties, Picard varieties and
N\'eron-Severi groups in the context of arithmetic geometry (i.e., when the base field is 
not assumed to be algebraically closed).  Apart from fixing notation, our goal here is
to record a technical fact about N\'eron-Severi groups of principal homogeneous spaces
(Proposition 18) which will come in handy later on.
\\ \\
For $V/k$ a variety, $\FAlb(V)$ denotes the total Albanese scheme
of $V$ and $\FPic(V)$ denotes the total Picard scheme of $V$.  These are the reduced
subschemes of the schemes parameterizing, respectively, zero-cycles on $V$ modulo rational
equivalence, and divisors modulo linear equivalence.  That is to say, these are objects
representing sheafified versions of the usual Albanese and Picard groups, so that 
e.g. $\FAlb(V)(k) = \Alb(V/\overline{k})^{\gk}$, and similarly for the Picard scheme.
One must keep in mind that the natural map $\Pic(V/k) \ra \FPic(V)(k)$ is injective
but not generally surjective (unless $V(k) \neq \emptyset$): that is, not every $k$-rational
divisor class need be represented by a $k$-rational divisor.  These two group schemes are locally algebraic; 
indeed, each is an extension of a finite rank
$\Z$-module by an abelian variety.  In the Albanese case this is induced by the
degree map:
\begin{equation}
0 \ra \FAlb^0(V) \ra \FAlb(V) \stackrel{\deg}{\ra} \Z \ra 0. 
\end{equation}
In the Picard case we have the subgroup scheme $\FPic^0(V)$ parameterizing divisor classes
algebraically equivalent to zero:
\begin{equation}
0 \ra \FPic^0(V) \ra \FPic(V) \ra NS(V) \ra 0;
\end{equation}
here $NS(V)$ is the N\'eron-Severi group of $V$.
\\ \\
We have that $\FAlb^0(V)/k$ and $\FPic^0(V)/k$ are abelian varieties that are
in duality.  Especially, if $A/k$ is an abelian variety, the 
map $P \mapsto [P] - [O]$ induces an isomorphism $A \stackrel{\sim}{\ra} \FAlb^0(A)$,
and this duality becomes the usual $\Pic^0(A) = A^{\vee}$.
\\ \\  
By takng $\overline{k}$-valued points in (1) or (2), we get short exact
sequences of $\gk$-modules.  The $\gk$-module structure on $\Z$ is
necessarily trivial, but it need not be so for the N\'eron-Severi group. 
As a matter of notation, we prefer to write $NS(V)$ for the Galois module
$NS(V)(\overline{k})$ and $NS(V)^{\gk}$ for $NS(V)(k)$.  
\\ \\
Now let $V/k$ be a principal homogeneous space for $A$, so $\FAlb^0(V) = A$,
$\FPic^0(V) = A^{\vee}$.  (In fact, a principal homogeneous space structure
on $V$ is equivalent to the choice of an isomorphism of $\FAlb^0(V)$ with
$A$.)  We have also that $V/\overline{k} \cong A/\overline{k}$, so that
$NS(A) \cong NS(V)$ as $\Z$-modules.  But more is true.
\begin{prop}
There exists a map $\psi: NS(A) \ra NS(V)$ which is an isomorphism of
$\gk$-modules.
\end{prop}
\noindent
Proof: Let $\mu: A \times V \ra V$ denote the $A$-action on $V$.
Also let $m: V \times V \ra A$ denote the corresponding subtraction map,
i.e., the $k$-map such that $m(\mu(a,v),v) = a$ for all $v \in V(\overline{k}), \
a \in A(\overline{k})$.  Fix any
$\overline{p} \in V(\overline{k})$, and consider
the isomorphism
\[\tau: V/\overline{k} \ra A/\overline{k}, \ v \mapsto m(v,\overline{p}). \]
We may use $\tau$ to pull back geometric line bundles from $A$ to $V$.  Since
the property of a line bundle being algebraically equivalent to zero is preserved
under all isomorphisms of abelian varieties, $\tau^*$ induces a $\Z$-module
isomorphism $\psi: NS(A) \ra NS(V)$.  We claim that $\psi$ is necessarily
$\gk$-equivariant, i.e., that for all $\sigma \in \gk$ and $L \in \Pic(A/\overline{k})$,
\[\psi(\sigma(L))-\sigma(\psi(L)) \in \Pic^0(V/\overline{k}). \]
Indeed
\[\psi(\sigma(L))-\sigma(\psi(L)) = 
\mu(\sigma(L),\overline{p})-\sigma(\mu(L,\overline{p})) = \]
\[\mu(\sigma(L),\overline{p})-\mu(\sigma(L),\sigma(\overline{p})) = \mu(\sigma(L),\overline{p})-
\mu((\mu(\sigma(L),\overline{p}),m(\sigma(\overline{p}),\overline{p})) \]
i.e., the difference between a line bundle and its translate, which is algebraically equivalent to zero.
\\ \\
Alternate proof: let $k(V)$ be the function field of $V$.  Then, since $V/k(V) \cong A/k(V)$, 
their N\'eron-Severi groups are isomorphic as $\ggg_{k(V)}$-modules, hence (since $k$ is algebraically closed
in $k(V)$) as $\gk$-modules.
\subsection{Polarizations versus strong polarizations}
Recall that a \textbf{polarization} on an abelian
variety $A/k$ is given by a geometric ample line bundle $L/\overline{k}$ which is algebraically equivalent
to each of its Galois conjugates: that is to say, the $k$-rationality condition on the polarization takes
place in $NS(A)$.  We must distinguish this from the notion of an element of $NS(A)$ represented by a 
$k$-rational ample line bundle, so we call the latter a \textbf{strong polarization}.  The obstruction to 
a polarization being strong comes from taking Galois cohomology of the exact sequence 
\[0 \ra A^{\vee}(\overline{k}) \ra \FPic(A)(\overline{k}) \ra NS(A) \ra 0 \]
to get
\[0 \ra A^{\vee}(k) \ra \Pic(A/k) \ra NS(A)^{\gk} \ra H^1(\gk,A^{\vee}(k)). \]
(Here we have used the fact that $\Pic(A/k) = \FPic(A)(k)$, since $A(k) \neq \emptyset$.)
In other words, to every $\lambda \in NS(A)^{\gk}$ we associate a class
$c_{\lambda} \in H^1(k,A^{\vee})$.  In fact, since $A^{\vee} = \Pic^0(A)$ classifies
skew-symmetric divisor classes on $A$, we have that $[-1]_A$ induces $-1$ on $\Pic^0(A)$,
whereas $[-1]_A$ acts as the identity on $NS(A)$, so that $c_{\lambda} \in H^1(k,A^{\vee})[2]$.
A thorough analysis of these classes is the subject of [Poonen-Stoll]; they show
in particular that $c_{\lambda}$ vanishes when $k$ is a $p$-adic field, so that when $k$ is a number
field $c_{\lambda} \in \Sha(k,A)[2]$ but (to say the least!) need not be zero in general.  
\\ \\
We will call an abelian variety $A/k$ \textbf{unobstructed} if $NS(A)^{\gk} \ra H^1(k,A^{\vee})$
is the zero map.
\subsection{Separate period-index problems for the Albanese and the Picard}
Let $V/k$ be any variety.  Define the \textbf{Albanese period} of $V$ to be the order
of the cokernel of the degree map $\FAlb(V)(k) \stackrel{\deg}{\ra} \Z$.  
For any $i \in \Z$, we denote $\FAlb^i(V)$ to be the degree $i$ component of the Albanese scheme, so $\FAlb^i(V)$ is 
a principal homogeneous space for $A: = \FAlb^0(V)$.  Then the Albanese period is the period of
the sequence $\FAlb^1(V), \ \FAlb^2(V) \ldots$, or the least positive $i$ such that
$\FAlb^i(V)(k) \neq \emptyset$.  In addition, using the exact sequence in Galois cohomology
derived from taking $\overline{k}$-valued points in (1), the Albansese period is
equal to the order of the kernel of the map $\Z \ra H^1(k,\FAlb^0(V))$.
\\ \\
Note that if $V(k) \neq \emptyset$ then the Albanese period of $V$ is equal to one, but not necessarily
conversely: e.g., (in characteristic zero) if $\pi_1(V/\overline{k}) = 0$, then the Albanese variety is
trivial, but $V$ need not have a rational point: one need look no further than genus zero curves.  On the other
hand, if $V/k$ is a principal homogeneous space for its Albanese variety, then the Albanese period is the
order of $V$ in $H^1(k,A)$ -- it is exactly the period in above cohomological sense.
\\ \\
The \textbf{Albanese index} of $V$ is the cokernel of the map $Z_0(V)^{\gk} \ra \Z$, the least positive
degree of a $k$-rational zero-cycle.  Notice that the least positive degree of an arbitrary zero-cycle
is the same as the greatest common divisor of degrees of \emph{effective} zero-cycles; it follows that
the Albanese index coincides with the index in the sense of Galois cohomology.
\\ \\
Define the \textbf{Picard period} of $V/k$ to be the \emph{exponent} of the cokernel of the map
$\FPic(V)(k) \ra NS(V)^{\gk}$.  Using the exact sequence (2), the Picard period is also the
exponent of the image of the connecting map $\delta: NS(V)^{\gk} \ra H^1(k,\FPic^0(V))$.  
Finally, define the \textbf{Picard index} of $V/k$ to be the exponent of the cokernel of
the map $\Div(V/\overline{k})^{\gk} \ra NS(V)^{\gk}$.  
\newcommand{\op}{\overline{p}}
\begin{prop}
Let $V/k$ be a principal homogeneous space of $A = \FAlb^0(V)/k$, and consider the map
$\delta: NS(V)^{\gk} \ra H^1(k,\FPic^0(V))$.  Suppose that $P \in NS(V)^{\gk}$ corresponds,
under the identification of $NS(V)$ with $NS(A)$, to the class of a $k$-rational line
bundle $L$ on $A$.  Then $\delta(P)$ is the image of $[V] \in H^1(k,\FAlb^0(V))$ under
the map $H^1(\varphi_{L}): H^1(k,\FAlb^0(V)) \ra H^1(k,\FPic^0(V))$.
\end{prop}
\noindent
Proof: We use the same notation as in the proof of Proposition 18; especially, recall
we have chosen a $\overline{p} \in V(\overline{k})$.  Our assumption is then
that $P$ is algebraically equivalent to a line bundle $L' = \psi(L) = \mu(L,\op)$
for some $L \in \Pic(A)$.  A cocycle representative for $\delta(P)$ is
then given by $\sigma \mapsto \sigma(L') - L'$, where
\[\sigma(L')-L' = \sigma(\mu(L,\overline{p}))-\mu(L,\overline{p}) = 
\mu(\sigma(L),\sigma(\overline{p}))-\mu(L,\overline{p}). \]
On the other hand, $\sigma \mapsto m(\sigma(\overline{p}),\overline{p})$ gives
a cocycle representative for  $V \in H^1(k,A)$.  The map
$\varphi_L: A \ra A^{\vee}$ is just $q \mapsto L_q - L$, so
the image of the cocycle under $H^1(k,A) \ra H^1(k,A^{\vee})$ carries $\sigma$ to
\[L_{m(\sigma(\op),\op)}-L. \]
Translating by $\overline{p}$ to view this as a line bundle on $V$, we get
$\mu(\sigma(L),\sigma(\op))-\mu(L,\op)$, completing the proof.
\begin{cor}
Let $V/k$ be a principal homogeneous space over an \emph{unobstructed} abelian variety $A/k$.
Then the Picard period of $V$ divides the Albanese period of $V$.  If $A/k$ is principally
polarizable then the Albanese and Picard periods of $V/k$ coincide.
\end{cor}
\noindent
\emph{Proof}: The proposition shows that we get a factorization
\[\delta: NS(V)^{\gk} \ra H^1(k,A)[n] \ra H^1(k,A^{\vee}) \]
so that the image of $\delta$ must be $n$-torsion.  Conversely,
a principal polarization induces an isomorphism $H^1(k,A) \ra H^1(k,A^{\vee})$
so that the image of $V$ has exact order $n$.  
\\ \\
Remark: The hypotheses can be weakened somewhat: it is enough to assume the
existence of a not necessarily positive $k$-rational line bundle $L$ such that $\varphi_L: A \ra A^{\vee}$ 
is an isomorphism.  Moreover, if the the Albanese period of $V$ is odd, one does not need to assume
that $A$ is unobstructed.
\\ \\
Remark: In any case, the proposition shows that if $V \in H^1(k,A)[n]$ and $P$ is a $k$-rational 
line bundle on $A$, then $n[P]$ is represented by a $k$-rational divisor \emph{class} on $V$.
\\ \\
Further relations betweeen the Picard period and the Picard index and between the Picard index and the 
Albanese index will be derived in the next section, using the period-index obstruction.  Putting these
together, we will get information on the relation between the Albanese period and the Albanese index,
which is what we really want.
\section{The period-index obstruction}
\noindent
 
\subsection{Definition of the period-index obstruction using theta groups}
The constructions of this section are the higher-dimensional analogues
of [O'Neil]. 
\\ \\
Let $A/k$ be an abelian variety, and let $L \in \Pic(A)$ be an ample line bundle, so that 
the associated map $\varphi_L : A \ra A^{\vee}$, $x \mapsto t_x^*(L) \otimes L^{-1}$ has
finite kernel, denoted $\kappa(L)$.  Let $\mathcal{G}(L)$ denote the theta group associated to $L$ 
(cf. [Mumford], [AV]);
it is the algebraic $k$-group representing the group functor of automorphisms of the total space of 
$L$ over translations of $A$.  The theta group lies in a short exact sequence
\[1 \ra \Gm \ra \mathcal{G}_L \ra \kappa(L) \ra 0. \]
Choose $f_1,\ldots,f_{N+1}$ a basis for the complete linear system associated to $L$, so that there
is an associated morphism $A \ra \mathbb{P}^N$.  Note that by Riemann-Roch, $N+1 = \sqrt{\#\kappa(L)}$. 
This choice of a basis allows us to identify
$\Aut(\mathbb{P}^N)$ with $PGL_{N+1}$; inside this full automorphism group we have the group 
of automorphisms preserving the image of $A$; one knows that this subgroup is precisely
$\kappa(L)$.  Accordingly, the theta group sequence maps to the short exact sequence
\[ 1 \ra \Gm \ra GL_{N+1} \ra PGL_{N+1} \ra 1 \]
the first vertical arrow being the identity, the last being the inclusion just described, and the
middle arrow being the map carrying $\alpha \in \Aut(L) \mapsto (\alpha^*f_1,\ldots,\alpha^*f_{N+1})$.
\\ \\
Enters Galois cohomology: having viewed $\kappa(L)$ as an automorphism group, we can now view
$H^1(k,\kappa(L))$ as a group of twisted forms.  Indeed, let $\phi: V \ra X$ be a $k$-morphism such
that $\phi/\overline{k}: A \ra \mathbb{P}^N$ is the morphism associated to $L$, up to an isomorphism
of $A$ extending to $\mathbb{P}^N$.  Then $H^1(k,\kappa(L))$ parameterizes these structures $\phi$.
Notice that $V$ is a twisted form of $A$; the assignment $(V \ra X) \mapsto V$ is seen on
cohomology as $H^1(k,\kappa(L)) \ra H^1(k,A)$.  On the other hand, $X$ is a twisted form of 
$\mathbb{P}^N$ -- i.e. a Severi-Brauer variety, and the assignment $(V \ra X) \mapsto X$ is
seen on cohomology as the connecting homomorphism $H^1(k,\kappa(L)) \stackrel{\delta}{\ra}
H^2(k,\Gm) = Br(k)$.  For $\xi \in H^1(k,\kappa(L))$ a cohomology class, we refer to the associated
element $\delta(\xi)$ of $Br(k)$ as the \textbf{period-index obstruction} associated to $\xi$.  
\\ \\
If $L' \in \Pic(A)$ is another line bundle which is algebraically equivalent to $L$ -- denoted
$[L] = [L']$, then $\kappa := \kappa(L) = \kappa(L')$ and we have two maps 
$\Delta_{L},\ \Delta_{L'}: H^1(k,\kappa) \ra Br(k)$.  Let us make explicit the relationship between these two maps.  
For this: Galois descent gives 
us a 
bijective correspondence between
$H^1(k,\kappa(L))$ and the set of equivalence classes of diagrams $\mathcal{D}_{L} = \{V \ra X\}$, 
where 
$V_1 \ra X_1$ is regarded as equivalent to $V_2 \ra X_2$ if there is an $A$-space isomorphism
$V_1 \ra V_2$ extending to an isomorphism of $X_1 \ra X_2$.  
\\ \\
Now observe that, as sets, $\mathcal{D}_{L} = \mathcal{D}_{L'}$ if and only if $[L] = [L']$; indeed
the ``only if'' is evident and the ``if'' can be read off from the Kummer sequence    
\[0 \ra A^{\vee}(k)/\phi_L(A(k)) \ra H^1(k,\kappa(L)) \ra H^1(k,A)[\phi_L] \ra 0. \]
This sequence also makes clear that the choice of an $L$ in its N\'eron-Severi class is
equivalent to the choice of an origin in $A^{\vee}(k)/\phi_L(A(k))$, and we find that
\[ \Delta_{L'}(\xi) = \Delta_{L}(\xi+[L'-L]). \]
In other words, the two obstruction maps differ only harmlessly -- especially, the question
of whether there exists a lift of $\eta \in H^1(k,A)[\phi]$ to $H^1(k,\kappa(L))$ with
vanishing obstruction is independent of the choice of $L$ in $[L]$.
\\ \\
Remark: The Kummer sequence also gives us the following useful fact: if $L$ is an ample line bundle on 
$A$ and $X$ is a principal homogeneous space over $A$ whose corresponding cohomology class is 
$\phi_L$-torsion, then there exists a rational divisor \emph{class} on $X$ representing the class 
$[L]$ in $NS(V) = NS(A)$; this was observed earlier by other means (see the remark following Corollary 21).
\subsection{Applications to the period-index problem}
\begin{thm}
Let $V/k$ be a principal homogeneous space of exact (Albanese) period $n$ over
a strongly principally polarized abelian variety $A/k$.  Suppose there exists
some lift of $V$ to $\xi \in H^1(k,A[n])$ such that $\Delta(\xi) = 0$ -- in particular
this occurs if the Picard index of $V$ is $n$.  Then $V$ can be split over
a field extension of degree at most $(g!) \cdot n^g$; in particular the
(Albanese) index of $V$ divides $g! \cdot n^g$.
\end{thm}
\noindent
\emph{Proof}: If $V/k$ is \emph{any} variety, there is a natural map of sets
\[\Theta: \mathcal{C} \ra \mathcal{P} \]
from the cone of ample $k$-rational divisors $\mathcal{C} \subset \Pic(V/k)$
to the subset $\mathcal{P} \subset \FAlb(V)(k)$ of $k$-rational points on 
the Albanese scheme which may be represented by effective zero-cycles.  Indeed,
we take the corresponding morphism $\varphi: V \ra \PP^N$ into projective space,
which has some degree $d$.  This means that if we cut with a
$k$-rational linear subspace whose dimension is equal to the codimension of
the image of $V$, we get an effective $k$-rational divisor of degree $d$ well-defined
up to rational equivalence, whence the map.  In our case, we may take
$D$ to be in the class of $nP$, so the degree of $\varphi$ is $(g!)\sqrt{\# \kappa(L)(\overline{k})} = (g!)\cdot n^g$.  
\\ \\
By assuming in addition the $\gk$-module triviality of $NS(A)$ (which, notice, is always satisfied in the
one-dimensional case), we can get a sort of converse result.
\begin{thm}
Let $(A,P)/k$ be a $g$-dimensional strongly principally polarized abelian variety and 
$V \in H^1(k,A)[n]$ a principal homogeneous space.  \emph{Suppose}
that the $\gk$-module action on $NS(A) = NS(V)$ is trivial.  If $\Delta(\xi) \neq 0$ for every lift
$\xi$ of $V$, then $V$ cannot be split over a degree $n$ field extension.  If $n= p^a$ is a prime power, we have
moreover that the index of $V$ exceeds $n$.
\end{thm}
\noindent
\emph{Proof}: We will prove the second statement first, so assume that $n = p^a$.  Seeking a contradiction, we 
assume that the index of $\eta$ is $p^a$, so that there
exists some splitting field $l_1/k$ of $\eta$ such that $[l_1:k] = m_1p^a$ with
$m_1$ prime to $p$.  On the other hand, since by Proposition 9 the index of $\eta$ is 
a power of $p$, for every prime $r$ dividing $m$ there exists a splitting extension $l_l/k$ such that 
$[l_l:k] = m_{r}$
for some $m_{r}$ prime to ${r}$.  Since $V/L_1 \cong A/L_1$, we may view the
strong principal polarization $P$ as an $L_1$-rational divisor on $V$.  
Let $D_1$ be the class of $Tr^{L_1}_k P$ in $\Pic(V/k)$.  By the assumed
$\gk$-module triviality of $NS(V)$, we must have $[D_1] = [m_1p^aP]$.
Repeating the previous sentence with $L_{r}$ in place of $L_1$, we get a
$D_2 \in \Pic(V/k)$ such that $[D_{r}] = [m_{r}P]$.  By varying over all ${r}$
dividing $m_1$, it follows that there exists $D \in \Pic(V/k)$ such that
$[D] = [p^aP]$, contradicting the nonvanishing of the period-index obstructions
of $\eta$.
\\ \\
A similar (but simpler) argument shows the first statement: if $l/k$ is a degree
$n$ field extension splitting $\eta$, then $Tr_{l/k} P$ exhibits a rational
divisor in the N\'eron-Severi class of $nP$.
\subsection{The proof of Theorem 3}
If $V \in H^1(k,A)[n]$ with $(A,P)$ a strongly principally polarized abelian variety, then in order to show that
$V$ can be split over a field extension of degree at most $(g!) \cdot n^g$ it is enough, by the preceding theorem,
to show that there exists some lift $\xi$ of $V$ with vanishing period-index obstruction.  Of course every lift
has this property if $Br(k) = 0$, so we look at the case where $k$ is a number field and $V$ has rational points
everywhere locally.  The proof is the same as in [O'Neil]: namely, at every completion $k_v$ of $k$,
$V(k_v) \neq \emptyset$ implies that the obstruction map $\Delta: \FPic(V)(k_v) \ra Br(k)$ vanishes identically.
(Alternately, the image under $\phi$ of a $k_v$-rational point of $V$ gives a $k_v$-rational point on the Severi-Brauer
variety $X$, making $X_V \cong \PP_v^N$.)  But the Hasse principle holds in the Brauer group of a number field, so
again \emph{any} lift of $V$ to $\xi \in H^1(k,A[n])$ has $\Delta_v(\xi) = 0$ for all $v$, hence
$\Delta(\xi) = 0$.
\\ \\
In fact this argument can be pushed further, using the reciprocity law in the Brauer group of a number
field.
\\ \\
Example (Kolyvagin classes): Let $\mathcal{K}(A/k)[n] \subset H^1(k,A)[n]$ be the subset of classes which are locally trivial
at all but at most one place of $k$.  Then, because the sum of the local invariants of a global Brauer
group element vanishes, the above proof still works to show that $\Delta$ vanishes identically
on the complete preimage of $\mathcal{K}(A/k)[n]$ in $H^1(k,A[n])$, so the mindex bound
of $(g!) \cdot n^g$ holds equally well for elements of $\mathcal{K}(A/k)[n]$.  In contrast to
$\Sha(A/k)[n]$, $\mathcal{K}(A/k)[n]$ can be infinite: notably, work of [Gross-Zagier] and [Kolyvagin] 
shows that $\mathcal{K}(E/k)[p]$ is infinite when $E/k$ is an elliptic curve over an imaginary quadratic
field of analytic rank $1$ and $p$ is a sufficiently large prime.  Recently, [Stein] has shown that there are genus one 
curves
over $\Q$ of every odd index; his construction can be seen as an exploitation of the fact that Kolyvagin's
classes necessarily have period equal to their index.
\\ \\

\section{Abstract theta groups and Galois cohomology}
\noindent
\subsection{}
An \textbf{abstract theta group} over a field $k$ is an algebraic $k$-group scheme $\mathcal{G}$ fitting into
a short exact sequence
\[1 \ra \Gm \ra \mathcal{G} \ra K \ra 0, \]
and satisfying the following additional properties: \\ \\
$\bullet$ The center of $\mathcal{G}$ is $\Gm$. \\
$\bullet$ $K$ is finite \'etale with underlying
abelian group of the form $A \oplus A$.  \\
$\bullet$ The characteristic of $k$ does not divide the order of $K$ (equivalently, of $A$).
\\ \\
Let $\delta = (d_1,\ldots,d_g)$ denote the sequence of elementary divisors
of the finite abelian group $A$; we say $\mathcal{G}$ is \textbf{of type}
$\delta$.  Let $d = \lcm(d_i)$ denote the exponent of $K$.  
If $d_1 = \ldots = d_g = n$, we will write $\underline{n}$ 
for $(n,n,\ldots,n)$.
\\ \\
We define, for any abstract theta group $\mathcal{G}$, the \textbf{commutator pairing} 
$e: K \times K \ra \Gm$, as follows: for any
$(P,P') \in K \times K$, lift to any elements $\tilde{P}, \tilde{P'}$
in $\mathcal{G}$.  Then $e(P,P')$ is defined to be the commutator
$[\tilde{P},\tilde{P'}]$, which is well-defined since the lift is ambiguous only
by central elements, and lives in $\Gm$ because of the commutativity of $K$.  It is
not hard to check that the condition that $\Gm$ be the precise center of $\mathcal{G}$
is equivalent to the \emph{nondegeneracy} of $e$.  So $e$ is a symplectic form,
taking values in $\mu_d \subset \Gm$, and placing $K$ into Cartier duality with itself.
\\ \\
For us, there are two important examples of abstract theta groups.  We have
already met the first: if $L$ is a nondegenerate line bundle on an abelian
variety $A/k$, then $\mathcal{G}_L$ is an abstract theta group.  Especially,
if $L = nP$ where $P$ is a strong principal polarization, then $\mathcal{G}_L$
has type $\underline{n}$.  
\\ \\
Here is another example: let $H/k$ be a finite \'etale group scheme of
order indivisible by the characteristic of $k$; put $K(H) = H \oplus H^*$
(remember that $H^* = \Hom(H,\Gm)$ denotes the Cartier dual of $H$).  There
is an abstract theta group denoted $\mathcal{H}(H)$, called the
\textbf{Heisenberg group associated to H}, defined as follows:
as a $\gk$-set, $\mathcal{H} = \Gm \times H \times H^*$, with the group law
given by
\[(\alpha,x,\ell) \star (\alpha',x',\ell') = (\alpha \alpha' \ell'(x),x+x',\ell + \ell'). \]
Note that in this
case the commutator pairing is $e((x,\ell),(x',\ell')) = x'(\ell)x(\ell')^{-1} \in \mu_d$.
We define the \textbf{standard Heisenberg group} $\mathcal{H}(\delta)$ of type $\delta$ to be the one
associated to the constant group scheme $\bigoplus_{i=1}^g \Z/d_i\Z$.  
\\ \\
Let us now explain why we have introduced this new terminology in the middle of our discussion
of the period-index osbstruction.  The point is that the cohomological obstruction map we have
defined for $\mathcal{G}_L$ makes sense for any abstract theta group, giving an ``abstract 
obstruction map''
\[\Delta_{\mathcal{G}}: H^1(k,K) \ra H^2(k,\Gm) = Br(k). \]
On one hand, we would certainly \emph{like} to see Heisenberg groups enter the picture, because in the Heisenberg
case the map $\Delta$ can be written down (more or less) explicitly: the most favorable
case is when $\mathcal{G} = \mathcal{H}(\underline{n})$ is the standard Heisenberg group and $k$ contains
the $n$th roots of unity, when (after an appropriate choice of basis) $\Delta$ is just a sum
of norm-residue symbols $\langle \ , \ \rangle_n$.  On the other hand, that Heisenberg groups
must play a role is guaranteed by the following result.
\begin{thm}([Mumford, Corollary to Theorem 1])
Suppose that $k = \overline{k}$ is algebraically closed.  Then every abstract theta group
of type $\delta$ is isomorphic to $\mathcal{H}(\delta)$.
\end{thm}
\noindent
Mumford's proof uses the surjectivity of multiplication by $d$ on $\Gm(k)$, so does not
go through for general $k$.  Indeed it is not true in general that an abstract
theta group needs to be isomorphic to a Heisenberg group, for two different reasons.
First, the definition of the Heisenberg group requires $K$ to 
be reducible as a $\Z/d\Z[\gk]$-module, but $A[n]$ is ``generically'' an irreducible
$\gk$-module.  More subtly, even when $A[n] \cong (\Z/n\Z)^{2g}$, it need not be the case
that $\mathcal{G}_{nP} \cong \mathcal{H}(\underline{n})$, as follows (e.g.) from the counterexample
discussed in [Clark].  We take instead a broader approach, using 
Galois cohomology to compare an \emph{abstract} theta group to a Heisenberg group (standard or
otherwise), and computing $\Delta_{\mathcal{G}}$ as a twisted form of the Heisenberg $\Delta$.
\subsection{The automorphism group}
Let $\Theta_{\delta}(k)$ be the
set of all abstract theta groups $\mathcal{G}$ of type $\delta$.  
By Mumford's theorem,
every element of $\Theta_{\delta}(k)$ is a $\overline{k}/k$-twisted form of $\mathcal{H}(\delta)$,
so $\Theta_{\delta}(k)$ should be identified as a Galois cohomology set $H^1(k,G_1)$,
where $G_1$ is a suitable automorphism group.  The goal of this section is to identify
$G_1$ as a subgroup of the full automorphism group of $\mathcal{H}(\delta)$ and to determine
its structure.
\\ \\
First we look carefully at the relevant descent problem: when we say that
$\mathcal{G}/\overline{k} \cong \mathcal{H}(\delta)/\overline{k}$, we mean by
an isomorphism $\iota$ which restricts to the identity map on $\Gm$.  Thus,
if $G_1 \subset \Aut(\mathcal{H})$ is group of automorphisms of the
Heisenberg group acting trivially on the center, we necessarily have
$\Theta_{\delta}(k) \subset H^1(k,G_1)$.  On the other hand, any
centrally trivial isomorphism $\iota: \mathcal{G} \ra \mathcal{H}$
induces an isomorphism $\overline{\iota}: K(\mathcal{G}) \ra K(\mathcal{H})$.  
Now $K(\mathcal{G})$ and $K(\mathcal{H})$ are both equipped
with symplectic forms $e$.  Since both symplectic forms are defined in terms
of the respective commutator pairings, a diagram chase reveals that 
$\overline{\iota}$ necessarily respects the symplectic structure.  It follows that 
every twisted form $\mathcal{T} \in H^1(k,G_1)$ is an abstract theta group,
and the map $H^1(k,G_1) \ra H^1(k,Sp(K))$ corresponds to
$\mathcal{G} \ra \mathcal{G}/\Gm$.
\\ \\
From now on we restrict our attention to type $\delta = \underline{n}$, the case
of interest to us in the sequel.  (The general case would only be notationally
more cumbersome.)
\\ \\
The next result gives a complete description of the group $G_1$.
\begin{prop}
Suppose that $\car(k) \neq 2$.  Then there is a split exact sequence
\[1 \ra K^* \ra G_1 \ra Sp(K) \ra 0. \]
\end{prop}
\noindent
Proof: The map $G_1 \ra Sp(K)$ is the one occurring in the definition of $G_1$; let $G_2$ be its kernel,
the group of automorphisms of $\mathcal{H}$ acting trivially on both $\Gm$ and the quotient $K$.  
\\ \\
Step 1: We claim that $G_2$ is canonically isomorphic to the character group of $K$.
Indeed, suppose that $\chi \in \Hom(K,\Gm)$ is any character of $K$.
Then $\chi: (\alpha,x,\ell) \mapsto (\chi(x,\ell)\alpha,x,\ell)$ gives an 
automorphism of $\mathcal{H}(\underline{n})$ acting trivially on
the center.  For the converse, since the Heisenberg group is
generated together by $\Gm$ together with any $\Z/n\Z$-basis $\{x_1,\ldots
x_g\} \cup \{\ell_1, \ldots, \ell_g\}$ for $H \oplus H'$, given an arbitrary $\psi \in G_2$, 
there is a unique character $\chi$ such that the action of $\chi^{-1} \circ \psi$ 
acts trivially on all elements of the form $(\alpha,x_i,0)$ and $(\alpha,0,\ell_j)$
hence is the identity map on all of $\mathcal{H}$.
\\ \\
Step 2: We will construct a section $Sp(K) \ra G_1$, which clearly suffices to prove the result.
The idea is as follows: the group law defining the Heisenberg group uses the bilinear form
$f((x,\ell),(x',\ell')) = \ell'(x)$, which is not such a nice form: it is degenerate and neither
symmetric nor alternating.  Moreover, we have in sight a symplectic form $e: K \times K \ra \Gm$.
So we wish we were dealing with the group $\mathcal{H'}$ constructed using $e$ instead of $f$, i.e.,
\[(\alpha,x,\ell) \star (\alpha',x',\ell') := (\alpha \alpha' e(\ell,\ell'),x+x',\ell + \ell'), \]
as in this case there is an evident faithful action of $Sp(K)$ by automorphisms
of $\mathcal{H}'$: \[g \mapsto ((\alpha,(x \oplus \ell)) \mapsto (\alpha,g(x \oplus \ell))). \]
But -- thanks to our assumption that $\car(k) \neq 2$ -- 
it turns out that $\mathcal{H} \cong \mathcal{H}'$.  Indeed,
let $W_1 = \{(1,x,0)\} \ W_2 = \{(1,0,\ell)\}$ be the standard level subgroups of $\mathcal{H}$,
so that every element of $\mathcal{H}$ has a unique expression of the form $w_1w_2\alpha$,
with $\alpha \in \Gm$.  We define the map
\[\Phi_{\jku}: w_1w_2\alpha \mapsto (\alpha + \frac{1}{2} e(\overline{w_1},\overline{w_2}),
\overline{w_1},\overline{w_2}). \]
Not only is $\Phi_{\jku}$ an isomorphism of groups preserving the $\gk$-module structures,
it is trivial on the center and on the quotient $K$, so that we may regard $G_1/G_2$ as being
\emph{canonically} isomorphic to $Sp(K)$.  This completes the proof.
\\ \\
Remark: The proof of Step 2 can be found in [Yu, $\S$ 10], except that Yu works 
with finite Heisenberg groups.  It is a key point for us that the same proof works for both cases; see Section 
6.5.
\subsection{The Lagrangian case}
We will now give an explicit computation of $\Delta$ when
$K(\mathcal{G})$ has a \textbf{Lagrangian Galois-module structure}.  
Namely, in this section we work with
an abstract theta group $\mathcal{G}$ which ``could be'' a Heisenberg
group in the sense that $\mathcal{G}/\Gm = K$ admits a $\gk$-module
decomposition $K = H \oplus H^*$ with $H$ and $H^*$ isotropic
subspaces for the symplectic form $e$.
\\ \\
Remark: In case $\mathcal{G} = \mathcal{G}_{p[O]}$ is the theta group
associated to a degree $p$ line bundle on an elliptic curve, this is equivalent
to assuming that the mod $p$ Galois representation is of \textbf{split Cartan} type.
When $k$ is a number field, such a Galois-module structure occurs for half of all
primes $p$ when $E$ has $k$-rational complex multiplication; otherwise, by Serre's theorem,
it can occur for at most finitely many $p$.  (Indeed, even for $E/K$ a Tate curve over a 
$p$-adic field, $E[n]$ can be Lagrangian for only finitely many $n$: see [Silverman,
Prop. V.6.1].)  For higher-dimensional abelian varieties, the assumption of a Lagrangian
structure on $A[n]$ is not quite so restrictive, but certainly one expects that for ``most'' abelian 
varieties over number fields and ``most'' positive integers $n$, $A[n]$ will be an irreducible $\gk$-module.
\\ \\
What is gained by restricting to the Lagrangian case is that $\mathcal{G}$ is a twisted form
of the corresponding Heisenberg group $\mathcal{H}(H)$ by an element of the smaller
group $H^1(k,G_2) = H^1(k,(H \oplus H^*)^*)$.  Let $\chi \in Z^1(k,G_2)$ be a one-cocycle
with values in the character group of $K(H)$.  We view $\mathcal{G}$ as a ``doubly twisted'' form
of $\Gm \times K(H)$: twisted as a group according to the cocycle $f$ as above, and with twisted
Galois-module structure using $\chi$.
\\ \\
We employ a more compact notation: write $P = (x,\ell) \in K(H)$, so that an arbitrary element
of $\mathcal{G}$ is written now as $(\alpha,P)$ and the group law is written as
$(\alpha,P) \star (\alpha',P') = (\alpha \alpha' f(P,P'),P+P')$, where $f(P,P') = \ell'(P)$ as before.
Note that $(\alpha,P)^{-1} = (\alpha^{-1}f(P,-P)^{-1},-P)$.
\\ \\
Now we can compute the coboundary map $\Delta_{\mathcal{G}}: H^1(k,K(H)) \ra H^2(k,\Gm)$ directly:
let $\eta \in Z^1(k,K(H))$.  Then
\[\Delta(\eta)(\sigma,\tau) = N_{\sigma} \sigma(N_{\tau}) N_{\sigma \tau}^{-1}, \]
where $N_{\sigma}, \ N_{\tau} , N_{\sigma \tau}$ are any lifts of $\eta(\sigma), \eta(\tau)$
and $\eta(\sigma \tau)$ to $\mathcal{G}$.  We choose the simplest possible lifts, namely
$N_{\sigma} = (1,\eta(\sigma))$ and so on.  So we get:
\[\Delta(\eta)(\sigma,\tau) = (1,\eta(\sigma)) \star \sigma(1,\eta(\tau)) \star (1,\eta(\sigma \tau))^{-1} = \]
\[(1,\eta(\sigma)) \star (\chi(\sigma)(\sigma(\eta(\tau))), \sigma(\eta(\tau))) \star
(f(\eta(\sigma \tau),-\eta(\sigma \tau))^{-1},-\eta(\sigma \tau)) = \]
\[(\chi(\sigma)(\sigma(\eta(\tau))) f(\eta(\sigma), \sigma(\eta(\tau))),\eta(\sigma) \sigma(\eta(\tau)) \star 
(f(\eta(\sigma \tau),-\eta(\sigma \tau))^{-1},-\eta(\sigma \tau)) =  \]
\begin{equation}
(\chi(\sigma)(\sigma(\eta(\tau))) f(\eta(\sigma),\sigma(\eta(\tau))),0). 
\end{equation}
Shortly we shall make further assumptions on the $\gk$-module structure of $K(H)$ that will allow us to write
this formula in a more useful way, but let us make two remarks in the present level of generality.
\\ \\
First, $\Delta(\eta)(\sigma,\tau)$ decomposes as a product of two terms $\Delta_1(\eta)\cdot \Delta_2(\eta)$.
The first term $\Delta_1$ is a \emph{linear form} in $\eta$, whereas the second term is a \emph{quadratic form}
in $\eta$, so altogether $\Delta$ is a quadratic map which is not necessarily a quadratic form -- it will be so
if and only if the character $\chi$ is trivial, i.e., if and only if $\mathcal{G} \cong \mathcal{H}(H)$ is a
Heisenberg group -- but it is certainly \emph{never} a linear form.  (This is not a surprise: [O'Neil] showed
that when $\mathcal{G}$ is the theta group of $nO$ on an elliptic curve, $\Delta$ is a quadratic map.)  
\\ \\
Second, both $\Delta_1$ and $\Delta_2$ visibly take values in the $n$-torsion subgroup of $Br(k)$, hence
so does $\Delta$ itself.  We isolate this for future reference:
\begin{prop}
Let $\mathcal{G}$ be a Lagrangian abstract theta group of type $\underline{n}$.  Then
\[\Delta_{\mathcal{G}}(H^1(k,K)) \subset Br(k)[n]. \]
\end{prop}
\subsection{The case of full level $n$ structure}
We now specialize further: we assume that $H \cong H^* \cong (\Z/n\Z)^g$ 
are \emph{both} constant: notice that this assumption implies that
$k$ contains the $n$th roots of unity.  So Kummer theory applies and we get
\[H^1(k,K) \stackrel{\beta}{\cong} H^1(k,\mu_n)^{2g} \cong (k^*/k^{*n})^{2g}.\]
The ``$\beta$'' denotes a choice of isomorphism, which is equivalent to a choice of
a $\Z/n\Z$-basis for $K(\overline{k})$.  In fact we want a careful choice of basis: first
fix $\zeta_n$ a primitive $n$th root of unity in $k$.  We choose a basis
$x_1,\ldots,x_d,y_1,\ldots,y_g$ which is ``$\zeta_n$-symplectic'' with respect to
$e$: for all $i$ and $j$, $e_n(x_i,x_j) = 0$, $e_n(x_i,y_j) = \delta_{ij} \zeta_n$. 
Then $\Delta = \Delta_{\mathcal{G}}$ may be viewed as a map
\[\Delta: (k^*/k^{*n})^{2g} \ra Br(k). \]
Let $(a_1,\ldots,a_d,b_1,\ldots,b_d)$ be an element of $(k^*/k^{*n})^{2g}$.
For $x,y \in k^*/k^{*n}$, let $\langle x, y \rangle_n \in Br(k)[n]$ denote the
norm-residue symbol [CL], and recall that the definition of the norm-residue symbol
requires a choice of a primitive $n$th root of unity (we choose the same $\zeta_n$).
\begin{thm}
For $1 \leq i \leq g$, there exist $C_{1,i}, \ C_{2,i} \in (k^*/k^{*n})$ such that
\[\Delta(a_1,\ldots,a_d,b_1,\ldots,b_d) = \sum_{i=1}^g \langle C_{1,i}a_i,C_{2,i} b_i \rangle
- \langle C_{1,i}, C_{2,i} \rangle =  \]
\[\left(\sum_{i=1}^g \langle a_i,b_i \rangle_n \right) + \left(\sum_{i=1}^g \langle C_{1,i},b_i \rangle_n +
\langle a_i,C_{2,i} \rangle_n \right). \]
\end{thm}
\noindent
Given what has already been said, the proof of this theorem is straightforward: the linear term $\Delta_1$
of equation (5) is a sum of $2g$ characters of $\gk$, so with Kummer-theoretic identifications is given
by elements $C_{1,i}, \ C_{2,i} \in k^*/k^{* n}$.  The quadratic term $\Delta_2$ is a sum of $g$ cup-products
of pairs of characters, so under the same identifications becomes the sum of $g$ norm-residue symbols (as usual).
For more details on the latter point -- and indeed, a complete proof that the Heisenberg $\Delta$ is a sum
of norm-residue symbols -- the reader may consult [Sharifi, Prop. 2.3].  Also the $\Delta_1$-part of the
computation is done in the one-dimensional case in [Clark].
\subsection{}
The aim of this section is to prove the following result, which is (almost) a generalization of Proposition 25.
\begin{thm}
Supppose that $\car(k) \neq 2$.  Let $n$ be any odd positive integer and $\mathcal{G}$ an abstract theta group of type $\underline{n}$.
Then $\Delta_{\mathcal{G}}(H^1(k,K)) \subset Br(k)[n]$.
\end{thm}
\noindent
Proof: By Section 6.2, we know that $\mathcal{G}$ corresponds to some twisted form of $\mathcal{H} = 
\mathcal{H}(\underline{n})$
via some cocycle $\varphi \in Z^1(k,G_1)$.  The idea is to show that, given the known structure of the finite
group $G_1$, twisting the coboundary map by $\varphi$ cannot twist us out of the $n$-torsion subgroup of $Br(K)$.
This is made rigorous as follows: we claim that there is a subgroup scheme $\mathcal{G}[n]$ of $\mathcal{G}$ whose
$\overline{k}$-valued points are the elements $(\alpha,x,\ell)$ with $\alpha \in \mu_n$.  Indeed, there is
certainly such a subgroup scheme of $\mathcal{H}(\underline{n})$, since $f(K,K) \subset \mu_n$.  It remains to
be checked that this subgroup of $\mathcal{H}(\overline{k}) = \mathcal{G}(\overline{k})$ is stable under
the $\varphi$-twisted $\gk$-action on $\mathcal{G}$, i.e., is stable under every possible automorphism 
$\varphi(\sigma) \in G_1$.  This can be restated as follows: every centrally-trivial automorphism of the 
Heisenberg group 
scheme with center $\Gm$ should restrict to an automorphism of the finite Heisenberg group scheme with
center $\mu_n$.  But $\mathcal{H}[n](\overline{k}) = H(n,g)$, a finite group of order $p^{2g+1}$, is exactly the object
studied in [Yu, $\S$ 10].  In Step 2 of the proof of Proposition 24, we noted that Yu's structure theory of
$H(n,g)$ carries over verbatim to $\mathcal{H}$: we now ask the reader to reverse this process  
and conclude that the group of centrally trivial automorphisms of $H(n,g)$ is $G_1$.  (Note that the ``$\frac{1}{2}$'' 
appearing in the argument now requires us to assume that $2$ is a unit in the ring $\Z/n\Z$.)  Thus the claim is true.
\\ \\
Thus there is such a sub-group scheme $\mathcal{G}[n]$ of $\mathcal{G}$.  Moreover, since $\mathcal{H}$ and
$\mathcal{H}[n]$ both have centrally-trivial automorphism group $G_1$, the cocycle $\varphi$ giving
$\mathcal{G}$ restricts to give the corresponding cocycle for $\mathcal{G}[n]$.  In other words, when
we go to compute the coboundary map $\Delta$, lifting by $P \mapsto (1,P)$ as in Section 6.3, then exactly the
same computation computes the coboundary map for the exact sequence
\[1 \ra \mu_n \ra \mathcal{G}[n] \ra K \ra 1. \]
since $\Delta_{\mathcal{G}[n]}$ visibly lands in $H^2(k,\mu_n) = Br(k)[n]$, the same must be true
of $\Delta_{\mathcal{G}}$, completing the proof.
\\ \\
Remark: The proof does not go through when $g = 1, \ n = 2$: the finite Heisenberg group $H(2,1)$ of order
$8$ can be shown to be isomorphic to the dihedral group $D_8$.  Since the center has order $2$,
the central-triviality condition is automatically satisfied, and $G_1$ for $H(2,1)$ is $\Aut(D_8) \cong D_8$.
So one gets that $G_1/G_2 \cong
\Z/2\Z$, whereas $Sp_2(\F_2) \cong S_3$.  
\\ \\
Nevertheless the theorem is still true in this case, in a stronger form: every Brauer group element in
the image of the obstruction map has \emph{index} $2$; see Section 6.7.  (But I don't know what happens when $n$ is even and
$g > 1$.)     
\subsection{Proof of Theorem 5}
Let $A/k$ be a strongly principally polarized abelian variety over a $p$-adic field, and $n$ a positive
integer.  If $n$ is even, we assume that $A[n]$ is a Lagrangian $\gk$-module.  By Proposition 25 and Theorem 27
we know that $\Delta(H^1(k,A[n])) \subset Br(k)[n]$.
\\ \\
Let $V \in H^1(k,A)[n]$ be any principal homogeneous space, and let $\xi$ be any lift
of $V$ to $H^1(k,A[n])$.  As above, $\Delta(\xi) \in Br(k)[n]$.  By hypothesis,
there exists a field extension $l/k$ of degree dividing $n^a$ such that
$\Delta(\xi)|_L = 0$.  Thus $V|L$ has trivial period-index obstruction, so by
Theorem 22, $V|L$ can be split over a field extension of degree at most
$(g!) \cdot n^g$.  Thus $V$ itself can be split over a field extension of
degree at most $(g!) \cdot n^{g+a}$.  
\subsection{Applications to the period-index problem in the Brauer group}
It is interesting to reconsider the examples of Lang and Tate (Proposition 12)
in light of Theorem 5.  Namely, let $k_g = \C((T_1))\cdots((T_{2g}))$,
let $(A,P)/\C$ be an abelian variety and view $A/k_g$ by basechange.  
The proof of Proposition 12 goes by 
observing that for any finite extension $l/k_g$, in the Kummer sequence
\[0 \ra A(l)/nA(l) \ra H^1(l,A[n]) \ra H^1(l,A)[n] \ra 0 \]
we have $A(l)/nA(l) = 0$; in other words, compatibly with 
restriction we have that $H^1(k_g,A)[n] \cong H^1(k_g,A[n]) \cong
(k_g^{*}/k_g^{*n})^{2g}$.  Thus, for any $1 \leq k \leq 2g$, the class
$(T_1,\ldots,T_k,1,\ldots,1)$ corresponds to an element $V$ of period
$n$ and index $n^{k}$.  Let $V$ be a class of index $n^{2g}$, $\xi$ the
(unique!) lift of $V$, and consider $\Delta(\xi)$.  By Propositon 24,
the period of $\Delta(\xi)$ divides $n$, whereas in Section 6.6
we saw that the index of $\Delta(\xi)$ divides $n^g$.  At least if
$n = p > g$, we \emph{must} have that the index of $\Delta(\xi) = n^g$,
since otherwise Theorem 5 would tell us that $\ord_p(i(V)) < p^{2g}$.
Thus we deduce period-index violations in the Brauer group of $k$.
Moreover, because $(A,P)$ arises
by basechange from an algebraically closed field, we have that 
$\mathcal{G}_{nP} \cong \mathcal{H}_g$,
so that in this case the period-index obstruction map is precisely the norm-residue
symbol.  It follows from the Merkurjev-Suslin theorem that, as $g \ra \infty$,
every class of period $n$ in $Br(k_g)$ appears as the period-index obstruction
of some principal homogeneous space $V/k_g$.
\section{The proof of Theorem 2}
\noindent
Let $(A,P)$ be a principally polarized abelian variety over a $p$-adic field
$k$ (recall that all abelian varieties over $p$-adic fields are unobstructed),
and let $V/k$ be a principal homogeneous space for $A/k$ of exact
(Albanese = Picard) period $n$.  We will show that the Picard index is $n$,
or more precisely that there exists a $k$-rational line bundle in the 
N\'eron-Severi class of $nP$ on $V/k$; as in Section 5, this gives us
a $k$-rational zero-cycle of degree at most $g!n^g$, which is what we want.
\\ \\
The true motivation behind the Albanese/Picard business of Section 3
can now be revealed: we want to be able to adapt the Cassels-Lichtenbaum
diagram to higher-dimensional principal homogeneous spaces.  Indeed, we have
the following, a direct generalization of [Lichtenbaum, p. 1216]: \\ \\
\newcommand{\Id}{\operatorname{Id}}
$\begin{CD}
Br(k)    @>>> Br(k) \\
@AAA                @AAA \\
0 \ra \FPic^0(V)(k) @>>> \FPic(V)(k) @>>> NS(A)^{\gk} @>>> H^1(k,\FPic^0(V)) \\
@AAA                @AAA @AAA \\
0 \ra \Pic^0(V) @>>> \Pic(V) @>>> NS(A)^{\gk}
\end{CD}$
\\ \\
\\ \\
Since $[P] \in NS(A)^{\gk} = NS(V)^{\gk}$, by Corollary 20 we have a $\gk$-invariant line 
bundle $D$ on $V/\overline{k}$ such that $[D] = [nP]$ in the N\'eron-Severi group.  
Suppose we can find a $D_0 \in \FPic^0(V)(k)$ such that $\Delta(D) = \Delta(D_0)$.
Then $\Delta(D-D_0) = 0$, $[D-D_0] = [D]$, and we would succeed in finding a
$k$-rational divisor in the N\'eron-Severi class of $nP$.
\\ \\
The key point is to know whether $\Delta(D) \in Br(k)[n]$.  If so, we finish by applying Theorem 3.ii
of [van Hamel], which says that the image of the obstruction map
$\Delta: \ \FPic^0(V)(k) \ra Br(k)$ has exact order $n$, so, using the fact
that $Br(k)[n] = \frac{1}{n}\Z/\Z$, we can find an invariant
line bundle algebraically equivalent to zero and with the same obstruction as
our $D$.  And indeed, by Proposition 25 and Theorem 27, we know that the image of the
obstruction map is contained in $Br(k)[n]$ under either of the hypotheses
of Theorem 2.  This completes the proof.
\section{Horizontal growth of the $p$-part of the Shafarevich-Tate group}
\noindent
In this section we give the proofs of Theorems 6 and 7 and of Corollary 8.  These results
were proved in [Clark] in the elliptic curve case.  
\\ \\
Suppose that $A/k$ is a strongly principally polarized abelian variety over a number
field and $p$ is a prime number such that $A[p]$ and $NS(A)$ are both trivial $\gk$-modules.
In order to show that a principal homogeneous space $V \in H^1(k,A)[p]$ has index exceeding $p$,
it suffices, by Theorem 23, to show that for every lift $\xi$ of $V$ to $H^1(k,A[n])$,
$\Delta(\xi) \neq 0$.  Using the Kummer sequence
\[0 \ra A(k)/pA(k) \ra H^1(k,A)[p] \ra H^1(k,A)[p] \ra 0\]
and the finiteness of $A(k)/pA(k)$ (weak Mordell-Weil theorem), the proof of Theorem 6 is reduced
to the following elementary claim.
\begin{prop}
Let $k$ be a number field containing the $p$th roots of unity, and let $H \subset (k^*/k^{* p})^{2g}$ be
any finite subgroup.  Then there exists an infinite subgroup $G \subseteq (k^*/k^{*p})^{2g}$ such that
for every nonidentity element $g$ of $G$ and every element $h \in H$, the Brauer group element
$\Delta(hg)$
is nonzero.
\end{prop}
\noindent
This is proven in [Clark] in the $g=1$ case, using standard results from local and global class field theory.
The general case can be proven in the same way -- in fact one can find such a subgroup consisting of
elements $(a_1,b_1,\ldots,a_d,b_d)$ for which $a_i = b_i = 1$ for $1 \leq i \leq d-1$, thereby reducing
to the one-dimensional case.
\\ \\
This completes the proof of Theorem 6.
\\ \\
For the proof of Theorem 7, let $G \subseteq H^1(k,A)[p]$ be an infinite subgroup each of whose nontrivial elements
cannot be split by a degree $p$ field extension.\footnote{Note that we showed the existence of an infinite subgroup
$G$ all of whose elements have index exceeding $p$, but we only use the weaker property that their mindex exceeds $p$.}  
Let $G' \subseteq G$ be a complementary subspace to the
finite-dimensional $\F_p$-subspace $\Sha(A/k)[p] \cap G$, so that $G'$ has finite index in $G$.  
We define a finite set $B \subset \Sigma_k$ of ``bad places''
of $k$ as follows: $B$ consists of all real places of $k$, all places lying over $p$ and all places at
which $A$ has bad reduction.  Define $H^1(k_B,A)[p] = \bigoplus_{v \in B}H^1(k_v,A)[p]$; this is a finite
abelian group.  So letting $k_{\Phi_1}$ be the kernel of the natural map $\Phi_1: H^1(k,A)[p] \ra H^1(k_B,A)[p]$,
we have that $G_1 := G' \cap k_{\Phi_1}$ has finite index in $G'$ and is therefore infinite.
\\ \\
Each nonzero element $\eta_1 \in G_1$ yields nontrivial elements of $\Sha$ over degree $p$
field extensions: let $\Sigma_1$ be the finite set of places of $k$ at which $\eta_1$ is locally
nontrivial, and consider any such place $v$.  By the theorem of Lang and Tate, $\eta_1 \in H^1(k_v,A)[p]$
can be split by a degree $p$ extension $l_v/k_v$, indeed by any ramified extension of degree $p$.  (More 
precisely,
we invoke [Lang-Tate] in the non-Archimedean case: if $v$ is
a real Archimedean place -- at which $g_1$ can only be nontrivial if $p = 2$ -- then obviously the class
splits over $\C$.)  By weak approximation, we can find infinitely many degree $p$ global extensions $l/k$ 
completing to each $L_v/k_v$ for all $v \in \Sigma_L$.  By construction, the restriction of $\eta_1$
to $l$ is everywhere locally trivial, so represents an element of $\Sha(A/L)[p]$.  Since $\eta_1$ does
not split over any degree $p$ field extension of $k$, $\eta_1$ represents a nonzero element in $\Sha(A/L)[p]$.
\\ \\
Finally, we can inductively build increasingly large finite subgroups $H_i$ of $G'$ such that the restriction
to suitable degree $p$ field extensions $l_i/k$ is injective on $H_i$, as follows: start with any
nonzero $\eta_1 \in G_1$ as above.  Put $B_2 := B \cup \Sigma_1$ and $G_2 := G_1 \cap k_{\Phi_2}$.  For
any nonzero element $\eta_2 \in G_2$, the set $\Sigma_2$ of places at which $\eta_2$ is locally nontrivial
is disjoint from $\Sigma_1$, so again we can find infinitely many degree $p$ global extensions with
prescribed global behavior at $\Sigma_1 \cup \Sigma_2$, and so on.  Continuing in this way, we can
for any $N$ construct a cardinality $N$ set of classes $\{\eta_1,\ldots,\eta_N\}$; since their sets
$\Sigma_1,\ldots,\Sigma_r$ are pairwise disjoint and nonempty, these classes automatically form an
$\F_p$-linearly independent set.  Let $l/k$ be any of the infinitely many global degree $p$ extensions
that simultaneously locally trivialize each $\eta_i$.  We finish by remarking that the $\eta_i$
remain $\F_p$-independent as elements of $\Sha(A/L)[p]$: if there exist $a_1,\ldots,a_N \in \F_p$ such that 
\[a_1\eta_1|_L + \ldots + a_N \eta_N|L = 0, \]
then the class $a_1\eta_1 + \ldots + a_N\eta_N$ is an element of $G$ which gets split over a degree $p$
field extension, so that $a_1 = \ldots = a_N = 0$.
\\ \\
Proof of Corollary 8: We proceed by extending the base field of an arbitrary principally polarized
abelian variety over a number field so that the hypotheses of Theorem 7 are satisfied.  That is,
given $A/k$, we need to bound the degree of a field extension $l_1/k$ such that 
namely, that $A[p]$ and $NS(A)$ are trivial Galois modules, and that the polarization comes from a
$l_1$-rational line bundle; we then apply Theorem 7 and pick up an extra factor of $p$.  First trivialize
$A[p]$ as a Galois module -- because of the Galois-equivariance of the Weil pairing,
the mod $p$ Galois representation lands in the general symplectic group $GSp_{2g}(\F_p)$.
The obstruction $c_{\lambda}$ to the principal polarization being a strong
polarization lives in $H^1(k,A)[2]$, so that by Proposition 9 this class can be split over an extension of
degree dividing $2^{2g}$.
It remains to trivialize the $\gk$-action on the N\'eron-Severi group.  We claim that the
$\gk$-action on $NS(A)$ of any $g$-dimensional abelian variety over any field $k$ can be trivialized
over a field extension of degree at most $\#GL_{4d^2}(\F_3)$.\footnote{At least in characteristic zero
this bound is certainly not optimal.  Our intent is merely to write down an explicit bound which
can be obtained without too much trouble.}  Indeed $NS(A) \cong \Z^g$, where
$g$ is the $\Q$ rank of Rosati-fixed subalgebra of $\End^0(A) = \End(A/\overline{k}) \otimes \Q$, so in any
case $g \leq \dim_{\Q} \End^0(A) \leq 4d^2$.  Letting $l/k$ be the splitting field for the Galois action,
we get a faithful representation $\rho: \ggg_{l/k} \hookrightarrow GL_{4d^2}(\Z)$, whence
$[L:k]$ is bounded above by the order of some finite subgroup of $GL_{4d^2}(\Z)$.  But the finite
subgroups of $GL_N(\Z)$ are uniformly bounded, and indeed -- see e.g. [LALG] -- there is the stronger fact that 
for any odd prime $\ell$, a finite subgroup $G$ of $GL_N(\Z_{\ell})$ has trivial intersection
with the kernel of reduction of $GL_N(\Z_{\ell}) \ra GL_N(\F_{\ell})$ hence has order at most $\#GL_N(\F_{\ell})$.
Taking $\ell = 3$ establishes the claim and completes the proof of Corollary 8.
\section{Final remarks}
\noindent
I. Combining Theorems 1 and 2, we get nearly optimal bounds on the period-index problem over $p$-adic fields: for
fixed $g$, and all odd primes $\ell > g$, the index divides the $g$th power of the period and when $\ell = p$
one can (for sufficiently large $p$-adic fields) do no better.  Looking however at the results of Lang-Tate
and Gerritzen, where a complete characterization of the splitting fields is known, one sees no additional trouble
coming from the ``small primes'' (not even from $2$), so it is natural to believe that the bound $i \ | \ n^g$
holds unrestrictedly.  One ought to look more closely at the case of good reduction and $n = p$, to which
formal group methods could presumably be applied.
\\ \\
II. Because of Theorem 5, for a good understanding of the period-index problem over more general fields $k$,
it would be very helpful to understand the period-index problem in the Brauer group of $k$.  Much work has been 
done on this, starting with the fundamental work of Hasse-Brauer-Noether on the case of a global field, and continuing
to the present day, with recent results due to [Saltman], [de Jong] and Lieblich.  In all cases, the results
obtained are of the form $i \ | \ n^{d-1}$, where $d$ is the cohomological dimension of $k$ in the sense of
[CG], and old (1935) examples due to Nakayama already show that indices as large as $n^{d-1}$ (for all $d$) arise
for suitable function fields.  A similar conjecture is the ``transitionality'' of the period-index problem: e.g..,
if $i \ | \ n^a$ for all classes in the Brauer group of $k$, then we expect $i \ | \ n^{a+1}$ for all classes in the
Brauer group of $k(t)$.  A general attack on either of these conjectures seems to be out of current reach.
\\ \\
III. We end by noting that the prime $2$ has played a distinguished role in several places in this paper: it came
up in the obstruction to a polarization being strong and again in studying automorphisms of Heisenberg groups.  Of
course $2$ is the most important prime in Mumford's theory of theta functions (whose surface we scuffed but did not
scratch in the present work).  And, for a final example, the recent paper [Polischuk] studies a property of abstract
theta groups which is automatically satisfied for both Heisenberg groups and theta groups of line bundles, namely the
existence of the Weil representation (given, in the theta group case, as the map $\mathcal{G}_L \ra GL_{n^g}$ of 
Section 5.1).  He shows that an abstract theta group need not admit a $k$-rational Weil representation, the obstruction
being a $2$-torsion element in $Br(k)$.  It would be very interesting to come up with a common framework explaining
each of these phenomena.

\end{document}